\documentclass[12pt]{article}

    \usepackage{amsmath}
    \usepackage{amsfonts}
\usepackage{color}
\usepackage{graphicx}

\newtheorem{theo}{Theorem}
\newtheorem{coro}{Corollary}

\newtheorem{prop}{Proposition}
\newtheorem{lemm}{Lemma}
\newtheorem{defn}{Definition}

\def\E{\mathbb{E}}
\def\S{\mathbb{S}}
\def\0{{\bf 0}}

\def\R{\mathbb{R}}

\def\cov{{\rm Cov}}

\def\B{{ B}}

\def\en2{{\epsilon_n^2}}
\renewcommand{\E}{\mathbb E \,}
\newcommand{\C}{{\cal C}}

\newcommand{\tod}{\stackrel{{\cal D}}{\longrightarrow}}
\newcommand{\eqd}{\stackrel{{\cal D}}{=}}

\newcommand{\eqco}{\setcounter{equation}{0}}
\newcommand{\thco}{\setcounter{theo}{0}}
\newcommand{\prco}{\setcounter{prop}{0}}
\newcommand{\laco}{\setcounter{lemm}{0}}
\newcommand{\coco}{\setcounter{coro}{0}}
\newcommand{\cjco}{\setcounter{conj}{0}}

\newcommand{\deco}{\setcounter{defn}{0}}
\newcommand{\allco}{\eqco  \thco \prco \laco \coco \cjco \deco}
\newcommand{\X}{{\cal X}}
\def\b{{\beta}}

\def\la{{\lambda}}

\def\v{{ w }}

\renewcommand{\P}{{{\cal P}}}

\newcommand{\Var}{{\rm Var}}

\newcommand{\Vol}{{\rm Vol}}

\newcommand{\liml}{\lim_{\lambda \to \infty} }

\newcommand{\tX}{{\tilde{X}}}
\newcommand{\tY}{{\tilde{Y}}}
\newcommand{\tQ}{{\tilde{Q}}}
\newcommand{\tP}{{\tilde{ {\cal P}}}}

\newcommand{\tp}{{\tilde{\phi}}}

\def\bdm{\begin{displaymath}}
\newcommand{\edm}{\end{displaymath}}
\def\benu{\begin{enumerate}}
\def\eenu{\end{enumerate}}
\def\beqn{\begin{equation}}
\def\eeqn{\end{equation}}
\def\be{\begin{equation}}
\def\ee{\end{equation}}
\def\bea{\begin{eqnarray}}
\def\eea{\end{eqnarray}}
\newcommand{\bean}{\begin{eqnarray*}}
\newcommand{\eean}{\end{eqnarray*}}
\newcommand{\bear}{\begin{eqnarray}}
\newcommand{\eear}{\end{eqnarray}}

\def\Comment#1{
\marginpar{$\bullet$\quad{\tiny #1}}}

\def\R{\mathbb{R}}
\def\B^2{\mathbb{D}}
\def\S{\mathbb{S}}
\def\B{\mathbb{B}}

\def\b{{\beta}}

\def\th{{\theta}}

\def\txi{{\xi^{(\la)} }}

\def\de{{\delta}}

\def\T{{T^{\la, z}}}

\def\Pl{{ {\cal P}^{\la, z}}}
\def\Plz{{ {\cal P}^{\la, z}_{r_z} }}
\def\pK{{ \partial K }}
\def\c{{ c^{\la, z} }}

\def\Comment#1{\lineskip-4pt
\marginpar{ $\bullet$\quad{\em\small #1}}}

\def\qed{\hfill\hbox{${\vcenter{\vbox{
    \hrule height 0.4pt\hbox{\vrule width 0.4pt height 6pt
    \kern5pt\vrule width 0.4pt}\hrule height 0.4pt}}}$}}

\def\la{{\lambda}}
\def\ka{{\kappa}}

\begin{document}
\title{\bf Variance Asymptotics and Scaling Limits for
Gaussian Polytopes}

\author{Pierre Calka$^{*}$,  J. E. Yukich$^{**}$ }

\maketitle

\footnotetext{{\em American Mathematical Society 2010 subject
classifications.} Primary 60F05, 52A20; Secondary 60D05, 52A23}

\footnotetext{$~^{*}$ Research partially supported by French ANR grant PRESAGE (ANR-11-BS02-003) and French research group GeoSto (CNRS-GDR3477)}\footnotetext{$~^{**}$ Research supported in part by NSF grant
DMS-1106619}

\begin{abstract} Let $K_n$ be the convex hull of i.i.d. random variables
distributed according to the standard normal distribution on $\R^d$.
We establish variance asymptotics as $n \to \infty$ for the re-scaled intrinsic volumes and
$k$-face functionals of $K_n$, $k \in \{0,1,...,d-1\}$, resolving an
open problem \cite{WW}.  Variance asymptotics are given in terms of functionals of germ-grain
models having parabolic grains with apices at a Poisson point
process on $\R^{d-1} \times \R$ with intensity $e^h dh dv$. The scaling limit of the boundary of $K_n$
as $n \to \infty$ converges
 to a festoon of parabolic surfaces, coinciding with that featuring in the geometric construction of the zero viscosity solution to Burgers' equation with random input. 
\end{abstract}

\section{Main results}\label{INTRO}

\allco

For all $\la \in [1, \infty)$, let $\P_\la$  denote a Poisson point process of intensity $\la
\phi(x)dx$, where $$\phi(x):= (2 \pi)^{-d/2} \exp(- {|x|^2 \over
2})$$ is the standard normal density on $\R^d, d \geq 2$. Let
$\X_n:= \{X_1,...,X_n \}$, where $X_i$ are i.i.d. with density
$\phi(\cdot)$.  Let $K_\la$ and $K_n$ be the Gaussian polytopes
defined by the convex hull of $\P_\la$ and $\X_n$, respectively.
The number of $k$-faces of $K_\la$ and $K_n$ are denoted by
$f_k(K_\la)$ and $f_k(K_n)$, respectively.

In $d = 2$, R\'enyi and Sulanke \cite{RS} determined $\E f_1(K_n)$ and later Raynaud
\cite{Ra} determined $\E f_{d-1}(K_n)$ for all dimensions.  Subsequently, work of
Affentranger and Schneider \cite{AS} and Baryshnikov and Vitale \cite{BV} yielded
the general formula
\be \label{GP}
\E f_k(K_n) = \frac{2^d} {\sqrt{d}} { \binom{d}{k + 1} }  \beta_{k, d - 1} ( \pi \log n)^{(d-1)/2}(1 + o(1)),\ee
with $k \in \{0,...,d-1\}$ and where  $\beta_{k, d - 1}$ is the
internal angle of a regular $(d-1)$-simplex at one of its
$k$-dimensional faces. Concerning the volume functional,
Affentranger \cite{Af} showed that its expectation asymptotics
satisfy \be \label{Aff} \E \Vol(K_n) = \kappa_d (2 \log n)^{d/2}(1 + o(1)),
\ee
where $\kappa_d:= \pi^{d/2}/\Gamma(1 + d/2)$ denotes the volume of the $d$-dimensional unit ball.

In a remarkable paper, B\'ar\'any and Vu \cite{BVu} use dependency
graph methods to establish rates of normal convergence for
$f_k(K_n)$ and $\Vol(K_n)$, $k \in \{0,...,d-1\}$.  A key part of
their work involves obtaining sharp lower bounds for  $\Var
f_k(K_n)$ and $\Var \Vol (K_n)$.  Their results stop short of determining
precise variance asymptotics  for $f_k(K_n)$ and $\Vol (K_n)$ as $n \to \infty$, an open
problem going  back to the 1993 survey of Weil and
Wieacker (p. 1431 of \cite{WW}).  We resolve this problem
in Theorems \ref{Th1} and \ref{Th2}, expressing variance asymptotics in terms
of scaling limit functionals of parabolic germ-grain models.

Let $\P$ be the Poisson point process on $\R^{d-1} \times \R$ with
intensity \be \label{defP} d\P((v,h)):= e^{h} dh dv, \ \  \text{with}
\ \ (v,h) \in \R^{d-1} \times \R.\ee
Let $\Pi^{\downarrow}:= \{(v,h) \in \R^{d-1} \times \R, h \le-
|v|^2/2\}$ and for $w:= (v,h) \in \R^{d-1} \times \R$ we
put $\Pi^{\downarrow}(w):= w \oplus
\Pi^{\downarrow}$, where $\oplus$ denotes Minkowski addition.  The maximal union
of  parabolic grains  $\Pi^{\downarrow}(w), w \in \R^{d-1} \times \R,$ whose interior  contains no point of $\P$ is
$$
\Phi(\P):= \bigcup_{\left\{ \substack{w\in \R^{d-1} \times\R \\
\P \cap \rm{int}(\Pi^{\downarrow}(w)) =\emptyset}\right.} \Pi^{\downarrow}(w).
$$
 Notice that
$\partial \Phi(\P)$ is a union of inverted parabolic surfaces.
Remove points of $\P$  not belonging to $\partial \Phi(\P)$
and call the resulting thinned point set ${\rm{Ext}}(\P)$.  See Figure 1.

We show that the re-scaled configuration of extreme points in $\P_\la$
(and in $\X_n$) converges to  ${\rm{Ext}}(\P)$ and that
 the scaling limit $\partial K_\la$ as $\la \to \infty$  (and of $\partial K_n$ as $n \to \infty$) coincides with  $\partial \Phi(\P)$.
Curiously, this boundary  features in the geometric construction
of the zero-viscosity solution of Burgers' equation \cite{Bur}.
We consequently obtain a closed form expression for
expectation and variance asymptotics for the number of shocks
in the solution of the inviscid  Burgers' equation,
adding to \cite{Ba}.

Fix $u_0:= (0,0,...,1)\in \R^d$
and let $T_{u_0}$ denote the tangent space to the unit sphere
$\S^{d-1}$ at $u_0$. The {\em exponential
 map} $\exp:= \exp_{d-1}: T_{u_0} \to \S^{d-1}$  maps a vector $v$ of 
 $T_{u_0}$ to the point $u \in \S^{d-1}$ such that $u$ lies at the end
 of the geodesic of length $|v|$ starting at $u_0$ and having initial direction $v.$
Here and elsewhere $| \cdot |$ denotes Euclidean norm.


\begin{figure}
\vspace{38mm}
\ \ \ \ \ \ \ \ \ \ \ \ \ \ \ \ \ \ \includegraphics[width=100mm]  {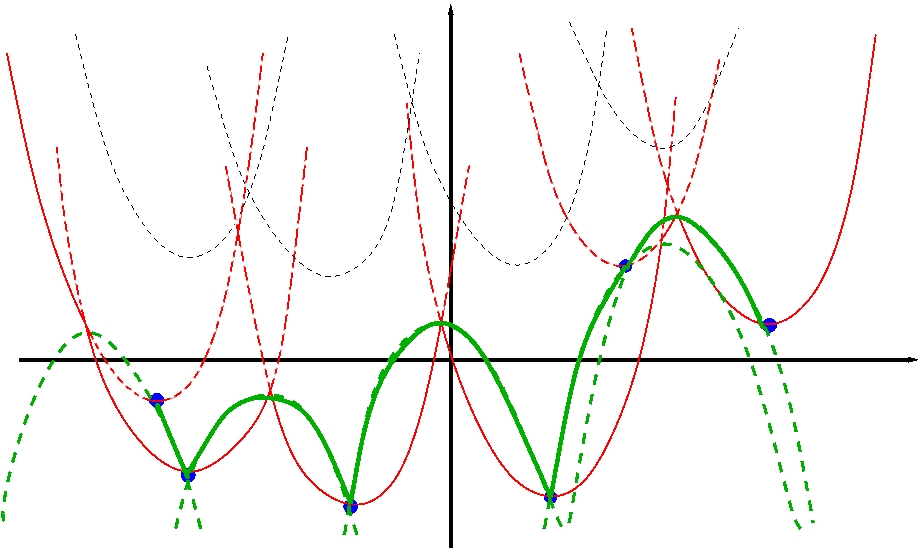}
\caption{The point process $\mbox{Ext}(\P)$ (blue);  the boundary of the germ-grain model $\partial(\Psi(\P)$ (red); the Burgers' festoon
 $\partial(\Phi(\P))$ (green). Points which are not extreme are at apices of gray parabolas. }
\end{figure}



For all $\la \in [1, \infty)$ put
\be \label{defR} R_\la := \sqrt{2 \log \la -
\log(2 \cdot (2 \pi)^d \cdot \log \la)}.\ee  Choose $\la_0$ so that for $\la \in [\la_0, \infty)$
we have $R_\la \in [1, \infty)$. Let $\B_{d-1}(\pi)$ be the closed Euclidean ball of radius $\pi$ and centered at the origin in the tangent space of $\S^{d-1}$ at the point $u_0$. It is also the closure of
the  injectivity region of $\exp_{d-1}$, i.e. $\exp(\B_{d-1}(\pi)) =
\S^{d-1}$. For  $\la \in [\la_0, \infty)$, define the scaling
transform
$T^{(\la)}: \R^d \to \R^{d-1} \times \R
$
by \be \label{scaltrans} T^{(\la)}(x) := \left(R_\la \exp_{d-1}^{-1}
{x \over |x| }, \ R_\la^2 (1 - {|x| \over R_\la})\right), \ \ x \in \R^d. \ee
Here $\exp^{-1}(\cdot)$ is the inverse exponential map, which is well defined
on $\S^{d-1} \setminus \{ - u_0 \}$ and which takes values in the injectivity region
$\B_{d-1}(\pi)$.  For formal completeness, on the `missing'
point $-u_0$ we let $\exp^{-1}$ admit an arbitrary value, say
$(0,0,\ldots,\pi),$ and likewise we put
$T^{(\la)}({\bf 0}) := (\0, R_\la^{2}),$
 where $\0$ denotes either the origin of $\R^{d-1}$ or $\R^d$,
 according to the context.

\noindent{Postponing the heuristics behind $T^{(\la)}$ until Section 3, we state our main results.}

\begin{theo} \label{Th0} Under the transformations $T^{(\la)}$ and $T^{(n)}$,  
the extreme points of the convex hull of the respective Gaussian samples $\P_\la$ and
$\X_n$ converge in distribution to the thinned process ${\rm{Ext}}(\P)$ as
$\la \to \infty$ (respectively, as $n \to \infty$).
\end{theo}

Let $B_d(v,r)$ be the closed $d$-dimensional Euclidean ball centered at $v \in \R^{d}$
and with radius $r \in (0, \infty)$.  $\C(B_{d}(v,r))$ is the space of
  continuous functions on $B_{d}(v,r)$ equipped with the supremum norm.

 \begin{theo}\label{Th4} Fix $L \in (0, \infty).$
 As $\la \to \infty$, the re-scaled boundary  $T^{(\la)}(\partial K_\la)$
 converges in probability to  $\partial (\Phi(\P))$
 in the space $\C(B_{d-1}(\0,L))$.
 \end{theo}
In a companion paper we shall show that
$\partial(\Phi(\P))$ is also  the scaling limit of the
boundary of the convex hull of i.i.d. points in polytopes. In $d = 2$, the reflection of $\partial(\Phi(\P))$ about the $x$-axis describes a
festoon of parabolic arcs featuring in the geometric construction of the zero viscosity solution($\mu = 0$) to Burgers' equation
\be \label{Bur}
{\partial v \over \partial t}  + (v, \nabla) v = \mu \Delta v,  \ \ v = v(t,x), \  t > 0, \ \ (x,t) \in \R^{d-1} \times \R^+,
\ee
subject to Gaussian initial conditions \cite{MSW}; see Remark (i) below.  Given its prominence in
 the asymptotics of Burgers' equation and its role in scaling limits of
 boundaries of  random polytopes,  we shall henceforth refer to $\partial (\Phi(\P))$ as  {\em the Burgers' festoon}.

The transformation $T^{(\la)}$ induces scaling
limit  $k$-face and volume functionals governing the large $\la$
behavior of convex hull functionals, as seen in the next results.
These scaling limit functionals are used in the description of the 
variance asymptotics for the $k$-face and volume functionals, $k \in \{0,1,...,d-1\}$.

\begin{theo} \label{Th1}
For all $k \in \{0,1,...,d-1\}$, there exists a constant $F_{k,d}
\in (0, \infty)$, defined in terms of averages of covariances of a scaling limit
$k$-face functional on $\P$,  such that
 \be \label{main2a}
\liml (2 \log \la)^{-(d-1)/2}  \Var f_k(K_\la) =
F_{k,d} \ee and
 \be \label{main2b}
\lim_{n \to \infty} ( 2 \log n)^{-(d-1)/2}  \Var f_k(K_n) = F_{k,d}.
\ee
\end{theo}

\begin{theo} \label{Th2}
There exists a constant $V_{d}
\in (0, \infty)$, defined in terms of averages of covariances of a scaling limit
volume functional on $\P$,  such that
 \be \label{maina}
\lim_{\la \to \infty} (2 \log \la)^{-(d-3)/2}  \Var \Vol(K_\la) =
V_{d} \ee and
 \be \label{mainb}
\lim_{n \to \infty} ( 2 \log n)^{-(d-3)/2}  \Var \Vol(K_n) = V_{d}.
\ee We also have 
 \be \label{mainexpect}
 \kappa_d^{-1} (2 \log \la)^{-d/2}  \E
\Vol(K_\la) = 1 -\frac{d\log(\log \la)}{4\log \la}+ O\left(\frac{1}{\log\la}\right).
  \ee
\end{theo}

The thinned point set ${\rm{Ext}}(\P)$ features in the description of asymptotic solutions to Burgers'
equation (cf. Remark (i) below) and we next consider its limit theory with respect to
the sequence of cylindrical windows
$Q_\la:= [-{1 \over 2} \la^{1/(d-1)}, {1 \over 2} \la^{1/(d-1)}]^{d-1} \times \R$ as $\la \to \infty$. 
The next result, a by-product of our general methods, 
 yields 
 variance and
expectation asymptotics for the number of points in ${\rm{Ext}}(\P)$ over
growing windows,  adding to \cite{Ba}.
\begin{coro} \label{Th3} 
There exist
constants $E_{d}$ and $N_{d} \in (0, \infty)$,
defined respectively in terms of averages of means and covariances of a thinning functional on  $\P$,
 such that
 \be \label{main1a}
\lim_{\la \to \infty}  
\la^{-1} \E [{\rm{card}}( {\rm{Ext}}(\P \cap Q_\la))] = E_{d} \ee and \be \label{main1b}
 \lim_{\la \to \infty}  
 \la^{-1}  \Var [{\rm{card}}({\rm{Ext}}(\P \cap Q_\la))] = N_{d}. \ee
 In particular, $N_{d} = F_{0,d}$.
\end{coro}

For $k\in\{1,\cdots,d-1\}$, we denote by $V_k(K_{\la})$ the $k$-th intrinsic volume of $K_{\la}$. In \cite{HR}, Hug and Reitzner establish  expectation asymptotics for $V_k(K_{\la})$ as well as an upper-bound for its variance.
The analog of  Theorem \ref{Th2} holds for $V_k(K_{\la})$, as shown by the next result, proved in Section \ref{sectionproofs}.
\begin{theo} \label{intrin} There exists a constant $v_k
\in [0, \infty)$, defined in terms of averages of covariances of a scaling limit
(intrinsic) volume functional on $\P$,  such that
\be \label{VarIV}
\lim_{\la \to \infty} (2 \log \la)^{-k+(d+3)/2}  \Var V_k(K_\la) = v_k.
\ee
Moreover,  we have
 $$
\frac{\kappa_{d-k}}{\binom{d}{k}\kappa_d} (2 \log \la)^{-k/2}  \E
V_k(K_\la) =
1 -\frac{k\log(\log \la)}{4\log \la}+ O\left(\frac{1}{\log\la}\right).
$$
\end{theo}
The limit \eqref{VarIV} improves upon Theorem 1.2 in \cite{HR} which shows that $(\log \la)^{-{(k-3)}/2}\Var V_k(K_{\la})$ is bounded. In \cite{ReBook}, Reitzner remarks `it seems that these upper bounds are not best possible'. 
We are unable to show that the limits $v_k$, $1\le k\le (d-1)$, are non-vanishing, that is to say we are unable to show optimality
of our bounds.  In particular, $\Var V_k(K_\la)$ goes to infinity for $k>(d+3)/2$ as soon as $v_k\ne 0$.

\vskip.5cm

There are several ways in which this paper differs from \cite{CSY}, which considers functionals of convex hulls on i.i.d. uniform points in $B_d(\0, 1)$.   First, as the extreme  points of a Gaussian sample are concentrated in the vicinity of the
critical sphere $\partial B_d(\0, R_\la)$,
 we need to calibrate the scaling transform $T^{(\la)}$ accordingly.  Second,
 the Gaussian sample $\P_\la$, when transformed by  $T^{(\la)}$, converges to a non-homogenous limit point process $\P$,
 which is carried by the whole of $\R^{d-1} \times \R$. This contrasts with \cite{CSY}, where the limit point process is
  simpler in that it is homogeneous and confined to the upper half-space. The non-uniformity of $\P$, together with its larger domain, induce spatial dependencies
 between the re-scaled functionals
 which are themselves non-uniform, at least  with respect to height coordinates.  The description of these dependencies is made explicit and  may be modified to describe the
 simpler dependencies of  \cite{CSY}. Non-uniformity of spatial dependencies leads to moment bounds for re-scaled $k$-face and volume functionals which are also non-uniform.
  Third,  the scaling limit of the boundary of the Gaussian sample converges to a festoon of parabolic surfaces, coinciding with that given by the geometric solution to Burgers' equation with random input.  This correspondence, described more precisely below, merits further investigation as it
suggests that some aspects of the convex hull geometry are captured by a stochastic partial differential equation.

\vskip.5cm

\noindent{\em Remarks.}
{\noindent}(i) {\em Burgers' equation.}  Let ${\rm{Ext}}(\P)'$
be the reflection of ${\rm{Ext}}(\P)$ about the hyperplane $\R^{d-1}$.
The point process  ${\rm{Ext}}(\P)'$ features in the solution to Burgers' equation
\eqref{Bur} for $\mu \in (0, \infty)$ as well as for $\mu = 0$ (inviscid limit).

When $\mu = 0$, $d = 2$,  and when the initial conditions are specified by
a stationary Gaussian process $\eta$ having covariance $\E \eta(\0) \eta(x) =
o(1/\log x), x \to \infty$, the re-scaled local maximum of the solutions converge in
distribution to ${\rm{Ext}}(\P)'$ \cite{MSW}.  The abscissas of points in  ${\rm{Ext}}(\P)'$ correspond to zeros of the limit velocity process $v(L^2t, L^2x)$, as $L \to \infty$
(here the initial condition is re-scaled in terms of $L$, not $L^2$).
See Figure 1 in \cite{MSW} as well as
Figure 13 in the seminal work of Burgers \cite{Bur}. The shocks in the limit velocity process coincide with the local
minima of the festoon $\partial(\Phi(\tilde{\P}))$, which are themselves the scaling limit of the projections of the origin onto the
hyperplanes containing the hyperfaces of $K_\la$.
By \eqref{scaltrans}, when $d=2$, the typical angular difference between consecutive extreme
points of $K_\la$, after scaling by $R_\la$, converges in probability to the
typical distance between  abscissas of points in  ${\rm{Ext}}(\P)'$.  {\em Thus the
re-scaled angular increments between consecutive extreme points in $K_\la$
behave like the spacings between zeros of the zero-viscosity solution to \eqref{Bur}.}

In the case $\mu \in (0, \infty)$, the point set ${\rm{Ext}}(\P)'$ is shown to be
the scaling limit as $t \to \infty$ of centered and re-scaled local maxima
of the solutions to Burgers' equation \eqref{Bur} when the initial conditions
are specified by degenerate shot noise with Poissonian spatial locations; see
Theorem 9 and Remark 3 of \cite{AMS}.   Correlation functions for ${\rm{Ext}}(\P)'$
are given in section 5 of \cite{AMS}.

\vskip.3cm


{\noindent}(ii) {\em Theorems \ref{Th0} and \ref{Th4} - related work.}
In 1961, Geffroy \cite{Ge} states that the Hausdorff distance between $K_n$ and
$B_d(\0, \sqrt{2\log n})$ converges almost surely to zero.  From \cite{BVu}  we also know  that the
extreme points of the polytope $K_\la$ concentrate around
the  sphere $R_\la \S^{d-1}$ with high probability. Theorems \ref{Th0}-\ref{Th4}
add to these results, showing convergence of the measure induced by the re-scaled extreme points as well as convergence of the re-scaled boundary.

\vskip.3cm

{\noindent}(iii) {\em Theorem \ref{Th1}- - related work.} As
mentioned, B\'ar\'any and \ Vu \cite{BVu} show \ that $(\Var
f_k(K_n))^{-1/2}(f_k(K_n) - \E f_k(K_n))$ converges to a normal
random variable as $n \to \infty$. They also show (Theorem 6.3 of
\cite{BVu})  that $\Var f_k(K_n) = \Omega(( \log n)^{(d-1)/2})$.
These bounds are sharp, as Hug and Reitzner \cite{HR} had previously
showed that $\Var f_k(K_n) = O(( \log n)^{(d-1)/2})$. 
Aside from
these variance bounds and  work of Hueter \cite{Hu}, asserting that
$\Var f_0(K_n) = c(\log n)^{(d-1)/2} + o(1)$, the second order
issues raised by Weil and Wieacker \cite{WW} have largely remained
unsettled in the case of Gaussian input.  In particular  the
question of showing
$$\Var f_k(K_n) = c(\log n)^{(d-1)/2}(1 + o(1))$$ for $k \in
\{1,...,d-1\}$ has remained open.  On page 298 of \cite{HR},  Hug
and Reitzner, commenting on the likelihood of progress, remarked
that `Most probably it is difficult to establish such a precise
limit relation...'. \  Theorem \ref{Th1} addresses these issues.


\vskip.3cm

{\noindent}(iv) {\em Theorem \ref{Th2}- -related work.} Hug and
Reitzner \cite{HR} show $\Var \Vol(K_n) = O(( \log n)^{(d-3)/2})$
and later B\'ar\'any and Vu \cite{BVu} show that $\Var \Vol(K_n) =
\Theta(( \log n)^{(d-3)/2})$.  
The asymptotics \eqref{maina} and \eqref{mainb} turn these bounds into precise limits.
The equivalence
\eqref{mainexpect} improves upon   \eqref{Aff} in the setting of Poisson input.



\vskip.3cm

{\noindent}(v) {\em Corollary  \ref{Th3}- -related work.}  Baryshnikov \cite{Ba} establishes  the asymptotic normality of $\text{card}({\rm{Ext}}(\P) \cap Q_\la)$ as $\la \to \infty$,
obtaining expectation and variance asymptotics in 
Theorem 1.9.2 of \cite{Ba}. Notice that ${\rm{Ext}}(\P) \cap Q_\la$ restricts extreme points in $\P$ to $Q_\la$, whereas ${\rm{Ext}}(\P \cap Q_\la)$
are the extreme points in $\P \cap Q_\la$, which in general is not the same set, by boundary effects.
 Baryshnikov left open the
question of obtaining explicit limits,
 remarking that `the
question of constants is quite tricky'; see p. 180 of ibid.

In general, if a point process $\P_{\infty}$ is  a scaling limit
to the solution of \eqref{Bur}, then  $\text{card} (\P_{\infty} \cap Q_\la)$ coincides with the number of Voronoi cells generated by the abscissas of points in $\P_{\infty} \cap Q_\la$;  under conditions on the viscosity and initial input, such cells model the matterless
voids in the Universe \cite{AMS, Ba, MSW}.

\vskip.3cm

{\noindent} (vi) {\em Goodman-Pollack model.} In view of the Goodman-Pollack model for Gaussian polytopes,  it is well-known
\cite{AS, BV,HR,ReBook} that asymptotics for functionals of $K_n$ admit counterparts for functionals
of the orthogonal projection of randomly rotated regular simplices in $\R^{n - 1}$.
The proof of \eqref{GP}, as given in \cite{AS}, is actually formulated as a limit result for the Goodman-Pollack model.
  Theorems \ref{Th1} and \ref{Th2}  may be likewise cast in terms of variances of projections
of high-dimensional random simplices. For more on the Goodman-Pollack model and its applications to coding theory, see  \cite{HR, ReBook}.

\vskip.3cm

This paper is organized as follows. Section \ref{Sec2} introduces scaling limit functionals of germ-grain models having parabolic grains.
These scaling limit functionals appear in a general theorem which extends and refines Theorems \ref{Th1} and \ref{Th2}.  In particular the limit
constants in Theorems \ref{Th1} and \ref{Th2} are seen to be the averages of scaling limit functionals on parabolic germ-grain models
carried by the infinite non-homogenous input $\P$.  Section \ref{Sec3} shows  for each $\la \in [1, \infty)$ that the scaling transform $T^{(\la)}$
maps the Euclidean convex hull geometry into `nearly' parabolic convex geometry, which in the limit $\la \to \infty$ becomes parabolic
convex geometry.  We show that the image of $\P_\la$ under $T^{(\la)}$ converges in distribution to $\P$ and that $T^{(\la)}$
defines re-scaled $k$-face and volume functionals.  Section \ref{Sec4} establishes that the re-scaled $k$-face and volume functionals
localize in space, which is crucial to showing the convergence of their means and covariances to the respective means and covariances of
their scaling limits. Finally Section \ref{sectionproofs} provides the proofs of the main results.




\vskip.3cm

\section{Parabolic germ-grain models and a general result} \label{Sec2}

\label{DefFestoon}
\allco

In this section we define scaling limit functionals of  germ-grain models
 and we use their second order correlations to precisely define the limit constants $F_{k,d}$ and $V_d$ in \eqref{main2a} and \eqref{maina}, respectively.
 We use the scaling limit functionals to
establish variance asymptotics for the empirical measures induced by the
$k$-face and volume functionals, thereby extending Theorems \ref{Th1} and \ref{Th2}. Denote points in $\R^{d-1} \times \R$ by $\v:= (v,h)$.

\vskip.5cm

\noindent{\bf 2.1. Parabolic germ-grain models}. Let $$\Pi^{\uparrow} := \{(v,h) \in \R^{d-1} \times \R^+, h \geq
{|v|^2 \over 2} \}.$$ Let $\Pi^{\uparrow}(w) := w \oplus \Pi^{\uparrow}$.
The point set $\P$  generates a germ-grain model of paraboloids
 $$
\Psi(\P):= \bigcup_{ \v \in \P } \Pi^{\uparrow}(\v).
$$
A point $w_0 \in \P$ is {\em extreme} with respect to $\Psi(\P)$ if
the grain $\Pi^{\uparrow}(w_0)$ is not a subset of the union of
the grains $\Pi^{\uparrow}(w), w \in \P \setminus \{w_0\}$. See Figure 1. It may be
verified that the extreme points from this construction coincide with ${\rm{Ext}}(\P)$, see e.g. section 3 of \cite{CSY}.

\vskip.5cm

\noindent{\bf 2.2. Empirical $k$-face  and volume
measures.} Given a finite point set $\X \subset \R^d$, let
$\rm{co}(\X)$ be its convex hull.

\begin{defn} \label{kface} Given $k \in \{0,1,...,d-1\}$
and $x$ a vertex of $\rm{co}(\X)$, define the $k$-face functional
$\xi_k(x, \X)$ to be the product of $(k +1)^{-1}$ and the number of
$k$-faces of $\rm{co}(\X)$ which contain $x$. Otherwise we put
$\xi_k(x, \X) = 0$. Thus the total number of $k$-faces in $\rm{co}(\X)$ is $\sum_{x \in \X}
\xi_k(x, \X)$. Letting $\de_x$ be the unit point mass at $x$, the empirical k-face measure for $\P_\la$ is \be
\label{zerom} \mu^{\xi_k}_\la:= \sum_{x \in \P_\la} \xi_k(x, \P_\la
) \de_x. \ee
\end{defn}

\vskip.5cm

Let ${\cal F}(x, \P_\la)$ be the collection of $(d-1)$-dimensional faces in $K_\la$ which contain
$x$ and let ${\rm{cone}}(x, \P_\la):= \{ry, r > 0, y \in {\cal F}(x, \P_\la) \}$ be the cone generated by ${\cal F}(x, \P_\la)$.

\begin{defn} \label{vol} Given  $x$ a vertex of \ $\rm{co}(\P_\la)$, define the
defect volume functional 
$$
\xi_V(x, \P_\la):= d^{-1}R_\la\left[ \Vol( {\rm{cone}}(x, \P_\la) \cap B_d(\0, R_\la)) -
\Vol( {\rm{cone}}(x, \P_\la) \cap K_\la ) \right].$$
When $x$ is not a vertex of \ $\rm{co}(\P_\la)$,  we put $\xi_V(x, \P_\la) = 0$.
The empirical defect volume measure is \be \label{volmeasure}
\mu^{\xi_V}_\la:= \sum_{x \in \P_\la} \xi_V(x, \P_\la) \de_x. \ee
\end{defn}
Thus the total defect volume of $K_\la$ with respect to the ball $B_d(\0, R_\la)$
is given by $R_\la^{-1}\sum_{x \in \P_\la} \xi_V(x, \P_\la).$

\vskip.5cm

\noindent{\bf 2.3. Scaling limit $k$-face and volume functionals}. A set of $(k+1)$ extreme
  points $\{w_1,...,w_{k + 1} \} \subset {\rm{Ext}}(\P)$,  generates a $k$-dimensional {\em parabolic face}
 of the Burgers' festoon  $\partial (\Phi(\P))$ if there exists a translate $\tilde{\Pi}^{\downarrow}$ of $\Pi^{\downarrow}$
such that $\{w_1,\cdots,w_{k+1} \}=  \tilde{\Pi}^{\downarrow} \cap
{\rm{Ext}}(\P)$. 
When $k = d -1$ the parabolic face is a
hyperface.

\begin{defn} \label{xiinf} 
Define the
scaling limit $k$-face functional $\xi^{(\infty)}_{k}(w, \P)$,  $k \in \{0, 1,..., d - 1\}$, to be the product of
$(k + 1)^{-1}$ and the number of $k$-dimensional parabolic faces of
the Burgers' festoon $\partial (\Phi(\P))$ which
contain $w$, if $w \in {\rm{Ext}}(\P)$ and zero otherwise.
\end{defn}

\begin{defn} \label{vol-xiinf} 
Define the
scaling limit defect volume functional $\xi^{(\infty)}_{V}(w,
\P), w \in {\rm{Ext}}(\P)$,  by
$$
\xi^{(\infty)}_{V}(w, \P):= d^{-1} \int_{ {\rm{Cyl}}(w)}
\partial (\Phi(\P))(v) dv,
$$
where  ${\rm{Cyl}}(w)$ denotes the projection onto  $\R^{d-1}$ of
the hyperfaces of $\partial( \Phi(\P))$ containing $w$.  Otherwise, when $w \notin {\rm{Ext}}(\P)$ we put
$\xi^{(\infty)}_{V}(w, \P)= 0$.
\end{defn}

One of the main features of our approach is that
$\xi^{(\infty)}_{k}, k \in \{0,1,...,d-1\},$ are scaling
limits of  re-scaled $k$-face functionals, as defined in Section 3.3.
A similar statement holds for $\xi^{(\infty)}_{V}$. Lemma \ref{L2} makes these assertions
precise.  Let $\Xi$ denote the collection of functionals
$\xi_{k}, k \in \{0,1,...,d-1\},$ together with
$\xi_{V}$.
Let $\Xi^{(\infty)}$ denote the collection of scaling limits
$\xi^{(\infty)}_{k}, k \in \{0,1,...,d-1\},$ together with
$\xi^{(\infty)}_{V}$.

\vskip.5cm

\noindent{\bf 2.4. Limit theory for  empirical $k$-face
and volume measures.} Define the following second order correlation
functions for $\xi^{(\infty)} \in \Xi^{(\infty)}$. 

\begin{defn} For all $\v_1, \v_2 \in \R^d$ and $\xi^{(\infty)} \in \Xi^{(\infty)}$ put
\be \label{SO2} c^{\xi^{(\infty)}}(\v_1,\v_2):=
c^{\xi^{(\infty)}}(\v_1, \v_2, \P):= \ee $$
 \E \xi^{(\infty)}(\v_1, \P
\cup \{\v_2\} ) \xi^{(\infty)}(\v_2, \P \cup \{\v_1\} ) -  \E
\xi^{(\infty)}(\v_1, \P ) \E \xi^{(\infty)}(\v_2, \P )$$
and  \be \label{S03}
 \sigma^2(\xi^{(\infty)}) := \int_{-\infty}^{\infty} \E \xi^{(\infty)}((\0,h_0), \P)^2 e^{h_0} dh_0  \ee
  $$ + \int_{-\infty}^{\infty}  \int_{\R^{d-1}} \int_{-\infty}^{\infty}
   c^{\xi^{(\infty)}}((\0,h_0),(v_1,h_1)) e^{h_0 + h_1} dh_1 dv_1 dh_0.
$$ 
\end{defn}
\vskip.5cm

Theorem \ref{Th1} is a special case of a general result
expressing the asymptotic behavior of the empirical $k$-face and volume
measures in terms of scaling limit functionals $\xi_k^{(\infty)}$
of parabolic germ-grain models. Let $\C(\S^{d-1})$ be the class of
bounded functions on $\R^d$ whose set of continuity points includes
$\S^{d-1}$. Given $g \in \C(\S^{d-1})$, let $g_r(x) := g(x/r)$ and let
$\langle g, \mu_\la^\xi \rangle$ denote the integral of $g$ with
respect to $\mu_\la^\xi$.
Let $\sigma_{d-1}$ be the $(d-1)$-dimensional surface measure on $\S^{d-1}$.
 The following is proved in Section \ref{sectionproofs}.

\vskip.5cm

\begin{theo} \label{Th5}  For all $\xi\in \Xi$
and $g \in \C(\S^{d-1})$   we
have \be \label{main1} \lim_{\la \to \infty} ( 2\log \la)^{-(d-1)/2}
\E [\langle g_{R_\la }, \mu_\la^{\xi} \rangle] =
\int_{-\infty}^\infty \E \xi^{(\infty)}((\0,h_0), \P) e^{h_0} dh_0
\int_{\S^{d-1}} g(u)
 d\sigma_{d-1}(u) \ee and \be \label{main2} \lim_{\la \to
\infty} ( 2 \log \la)^{-(d-1)/2}  \Var [\langle g_{ R_\la },
\mu_\la^{\xi} \rangle] = \sigma^2(\xi^{(\infty)})
\int_{\S^{d-1}} g(u)^2
 d\sigma_{d-1}(u) \in (0, \infty). \ee
\end{theo}

\vskip.5cm





\vskip.5cm

{\noindent}{\em Remarks.}
{\noindent}(i) {\em Deducing Theorems \ref{Th1} and \ref{Th2} from
Theorem \ref{Th5}}. Setting $\xi$ to be $\xi_k$, the convergence \eqref{main2a}
is implied by \eqref{main2} with
$F_{k,d}=\sigma^2(\xi^{(\infty)}_{k}) \cdot d\ka_d$, with $d\ka_d = d \pi^{d/2}/\Gamma(1 + d/2)$ being the surface area of the unit sphere. Indeed, applying
\eqref{main2} to $g \equiv 1$, we have
$$\langle 1, \mu_\la^{\xi_k} \rangle=\sum_{x\in \P_{\la} }\xi_k(x, \P_\la
) =f_k(K_{\la}).$$
Likewise, putting $g \equiv 1$ in \eqref{main2}, setting $\xi$ to be $\xi_V$,
and recalling that $\xi_V$ incorporates an extra factor of $R_\la$, we get
the convergence \eqref{maina}, with $V_d := \sigma^2(\xi^{(\infty)}_{V})\cdot d\ka_d$.

To obtain \eqref{mainexpect},  put $g \equiv 1$ in \eqref{main1} and set $\xi \equiv \xi_V$ to get
$(2 \log \la)^{-d/2} \E [ \Vol (B_d(\0, R_\la)) - \Vol(K_\la) ] = O(R_\la^{-1} (\log \la)^{-1/2} )= O((\log \la)^{-1})$.
We have $$R_\la=\sqrt{2\log \la}-\frac{  \sqrt{2} \log(\log \la)}{4\sqrt{\log \la}}+ O \left(\frac{1}{\sqrt{\log\la}}\right)$$
which gives $(2 \log \la)^{-d/2} R_\la^d = 1 -d\log \log \la/ (4\log \la)+O((\log\la)^{-1})$ and thus
\eqref{mainexpect} holds.
When $\xi$ is set to $\xi_V$, we are unable to show that the right side of \eqref{main1} is non-zero,
that is  we are unable to show
$(2 \log \la)^{-d/2} \E [ \Vol (B_d(\0, R_\la)) - \Vol(K_\la)] = \Omega(( \log \la)^{-1})$.

The de-Poissonized limit \eqref{mainb} follows from the coupling of binomial
and Poisson points used in B\'ar\'any and Vu \cite{BVu}, in particular Lemma 8.1 of
\cite{BVu}.  The limit  \eqref{main2b} similarly follows from \eqref{main2a} and the same coupling, as described
in Section 13.2 of \cite{BVu}.


\vskip.3cm

{\noindent}(ii) {\em Central limit theorems.}  Combining \eqref{main2} with the
results of \cite{BVu} shows the following central limit theorem, as
$\la \to \infty$:
\be \label{CLT1}
( 2 \log \la)^{-(d-1)/2}(\langle g_{ R_\la }, \mu_\la^{\xi_k}
\rangle - \E \langle g_{ R_\la }, \mu_\la^{\xi_k} \rangle) \tod
N(0, \sigma^2),
\ee
where $N(0, \sigma^2)$ denotes a mean zero normal random variable
with variance  $\sigma^2:=
\sigma^2(\xi^{(\infty)}_{k}) \int_{\S^{d-1}} g(u)^2 du.$
Alternatively, using the localization of functionals  $\xi \in \Xi$, as described in Section 4,
together with standard stabilization methods as in \cite{CSY}, we obtain another proof of \eqref{CLT1}.

\vskip.5cm

\noindent{\bf 2.5. Further extensions.}
\noindent{\em (i) Brownian limits.}
Following the scaling methods
  of this paper and by appealing to the methods of section 8 of \cite{CSY} we may deduce that
the process given as the integrated version of the defect volume converges to
a Brownian sheet process.  This goes as follows.
  For $\X \subset \R^d$ and $u \in \S^{d-1}$ we put
$$
r(u, \X):= R_\la - \sup\{ \rho > 0: \rho u \in  {\rm{co}}(\X) \} 
$$
and for all $\la \in [\la_0, \infty)$ let $r_\la(u) := r(u, \P_\la)
$. 
Recall that $\B_{d-1}(\pi)$ is the closure of
the  injectivity region of $\exp_{d-1}$.
Define  for $v \in \B_{d-1}(\pi)$  the
{\em defect volume process}
$$
V_\la(v):= \int_{ \exp([\0, v])} r_\la(u) d \sigma_{d-1}(u).
$$
Here   $[\0,v]$ for $v \in \R^{d-1}$ is the rectangular solid
 in $\R^{d-1}$ with vertices $\0$ and $v,$ that is to say $[\0,v] :=
 \prod_{i=1}^{d-1} [\min(0,v^{(i)}),\max(0,v^{(i)})]$, with $v^{(i)}$ standing
 for the $i$th coordinate of $v$.  When $[T^{(\la)}]^{-1}[\0,v] = \B_{d-1}(\pi)$, we have that
 $V_\la(v)$ is the total defect volume of $K_\la$ with respect to $B_d(\0, R_\la)$.
We re-scale $V_\la(v)$ by its standard deviation, which in view of \eqref{maina}, gives
$$
\hat{V}_\la(v):= (2 \log \la)^{-(d - 3)/4 }(V_\la(v) - \E V_\la(v)), \ v \in \R^{d-1}.
$$
For any $\sigma^2 > 0$ let $B^{\sigma^2}( \cdot)$ be the Brownian sheet of variance
 coefficient $\sigma^2$ on the injectivity region $\B_{d-1}(\pi)$.
Extend the domain of  $B^{\sigma^2}$ to all of  $\R^{d-1}$ by
defining $B^{\sigma^2}$ as
the mean zero continuous path
 Gaussian process indexed by $\R^{d-1}$ with
 $$ \cov(B^{\sigma^2}(v),B^{\sigma^2}(w)) = \sigma^2 \cdot \sigma_{d-1}(\exp([\0,v] \cap [\0,w])).
 $$

  \begin{theo}\label{brown}
  As $\la\to\infty$, the random functions $\hat{V}_{\la} : \R^{d-1} \to \R$ converge in law to the Brownian sheet $B^{\sigma^2_V}$ in ${\cal C}(\R^{d-1}),$
  where $\sigma^2_V:= \sigma^2(\xi_V^{(\infty)})$.
 \end{theo}

We shall not prove this result, as it follows closely the proof of
Theorem 8.1 of \cite{CSY}.

\vskip.5cm

\noindent{\em (ii) Binomial input.} By coupling binomial
and Poisson points as in \cite{BVu}, we deduce the binomial analog
of Theorem \ref{Th5} for  measures $\sum_{i=1}^n
\xi(X_i, \X_n) \delta_{X_i}, \xi\in \Xi$,  where we recall that
$X_i$ are i.i.d. with density $\phi$ and $\X_n := \{X_j \}_{j=1}^n.$

\vskip.5cm

\noindent{\em (iii) Random polytopes on general Poisson input.} We expect that our main results extend to random polytopes generated by Poisson points
having an isotropic  intensity density.  As shown by Carnal \cite{Ca} and others,
there are qualitative differences in the behavior of $\E f_k(K_n)$ according to
whether the input $\X_n$ has an exponential tail or an algebraic tail modulated by a slowly varying function.
The choice for the critical radius $R_\la$ and the scaling transform $T^{(\la)}$ would thus need to reflect such
behavior.  For example, if $X_i, i \geq 1$, are i.i.d. on $\R^d$ with an isotropic intensity density decaying exponentially
with the distance to the origin and if
 $R_\la = \log \la - \log \log \la$,
then $T^{(\la)}(\P_\la) \tod {\mathcal H}_1$, where ${\mathcal H}_1$ is a rate one homogenous Poisson point process on $\R^d$.

\section{Scaling transformations} \label{Sec3}

\allco

For all $\la \in [\la_0, \infty)$, the scaling transform $T^{(\la)}$ defined at \eqref{scaltrans} maps $\R^d$ onto the
rectangular solid $W_\la \subset \R^{d-1} \times \R$ given by
$$W_\la:= (R_\la \cdot \B_{d-1}(\pi))
\times (-\infty, R_\la^2].$$ 
Let $(v,h)$ be the coordinates in $W_\la$, that is
\be \label{defvh}
v = R_\la \exp_{d-1}^{-1} { x \over |x|}, \ \ h = R_\la^2(1 - {|x| \over R_\la}).\ee
 Note that
$\S^{d-1}$ is geodesically complete in that $\exp_{d-1}$ is well defined on the whole tangent space
 $\R^{d-1} \simeq T_{u_0},$ although it is injective only
 on $\{ v \in T_{u_0},\; |v| < \pi \}.$

 The reader may wonder about the genesis of $T^{(\la)}$ and the parabolic scaling by $R_\la$.
Roughly speaking,
the effect of   $T^{(\la)}$ is to first re-scale the Gaussian sample by the characteristic scale factor
$R_\la^{-1}$ so that $\partial K_\la$ is close to $\S^{d-1}$. By considering the distribution of
$\max_{i \leq n} |X_i|$ we see that $(1 - |x|/R_\la)$ is small when $x \in \partial K_\la$; cf. \cite{Ge}. Re-scale again according to the twin desiderata:
(i) unit volume image subsets near the hyperplane $\R^{d-1}$ should host $\Theta(1)$ re-scaled points, and
(ii) radial components of points should scale as the square of angular components
$\exp_{d-1}^{-1}{x/ |x|}$.
Desideratum (ii) preserves the
parabolic behavior of the defect support function for $R_\la^{-1} K_\la$,
namely the function $1 - h_{R_\la^{-1} K_\la}(u), u \in \S^{d -1}$, where $h_K$ is the support function of $K \subset \R^d$.  Extreme
value theory  \cite{Res} for $|X_i|, \ i \geq 1$, suggests (i) is achieved via radial scaling by $R_\la^2$, whence by (ii) we obtain angular scaling of $R_\la$, and \eqref{scaltrans} follows.  These heuristics are justified below, particularly through Lemma \ref{weakcon}.
In this and in the following section, our aim is to show:
\vskip.3cm

{\noindent}(i) $T^{(\la)}$ defines a $1-1$ correspondence between
boundaries of convex hulls of point sets $\X \subset \R^d$ and a subset of
piecewise smooth functions on $W_\la,$

\vskip.1cm

{\noindent}(ii) 
$T^{(\la)}(\P_\la)$ converges in distribution to $\P$ defined at
\eqref{defP}, and

\vskip.1cm

{\noindent}(iii) $T^{(\la)}$ defines re-scaled $k$-face and volume
functionals on input carried by $W_\la$; when the input is $T^{(\la)}(\P_\la)$ then as $\la \to \infty$ the
means and covariances converge to the respective means and
covariances of the corresponding functionals in $\Xi^{(\infty)}.$

\vskip.5cm

\noindent{\bf 3.1. The re-scaled boundary of the convex hull under
$T^{(\la)}$}. 
Abusing notation, we let $\langle \cdot, \cdot \rangle$ denote inner product on
$\R^d$. For $x_0 \in \R^d \setminus \{\0 \},$ consider the ball
$$B_d({x_0 \over 2}, {|x_0| \over 2} )=\{x\in\R^d\setminus\{\0\}: \ |x| \leq \langle x_0, \frac{x} {|x|} \rangle \}\cup\{\0\}.$$
Consideration of the support function of $\rm{co}(\X)$ shows that $x_0 \in \X$ is a vertex of $\rm{co}(\X)$
iff   $B_d(x_0/2, |x_0|/2)$  is not a subset of $\bigcup_{x \neq x_0} B_d(x/2, |x|/2)$.
With $d_{\S^{d-1}}$ standing for the geodesic distance in
$\S^{d-1}$, let $\theta:= d_{\S^{d-1}} ( {x / |x|}, {x_0 / x_0|} )$. We  rewrite $B_d(x_0/2,|x_0|/2)$ as
$$
B_d({x_0 \over 2},{|x_0| \over 2})= \{x \in \R^d: \ |x| \leq |x_0| \cos \theta \}.
$$
Recalling the change of variable at \eqref{defvh}, let
$T^{(\la)}(x_0) := (v_0, h_0)$, so that
$h_0 = R_\la^2(1 - |x_0|/R_\la).$
We may rewrite  $B_d(x_0/2,|x_0|/2)$ as
$$
B_d({x_0 \over 2}, {|x_0| \over 2}):= \{x \in \R^d: \ R_\la^2(1 - {|x| \over R_\la \cos \theta}) \geq R_\la^2(1 - {|x_0| \over R_\la}) \}.
$$
Thus for all $\la \in [\la_0, \infty)$, $T^{(\la)}$ transforms $B_d(x_0/2,|x_0|/2)$ into
  the upward opening grain
\be \label{formu4.2}
[\Pi^{\uparrow}(v_0, h_0)]^{(\la)} := \{ (v,h) \in W_\la, \;
     h \geq  R_\la^2 (1-\cos[e_\la(v, v_0)]) + h_0\cos[e_\la(v,v_0)] \},
     \ee
with \be \label{defe} e_\la(v, v_0):=
d_{\S^{d-1}}(\exp_{d-1}(R_\la^{-1} v), \exp_{d-1}(R_\la^{-1} v_0)).
\ee 
Every finite $\X \subset W_\la$, $\la \in [\la_0, \infty)$, generates
the germ-grain model
\be \label{Boolmod}
\Psi^{(\la) }(\X):= \bigcup_{w \in \X } [\Pi^{\uparrow}(w)]^{(\la)}.
\ee
 This germ-grain model has a twofold relevance: (i)
$[\Pi^{\uparrow}(T^{(\la)}(x)) ]^{(\la)}$ is the image by $T^{(\la)}$ of $B_d(x/2, |x|/2)$ and (ii) $x \in \X$ is a vertex of $\rm{co}(\X)$ if and only
if $[\Pi^{\uparrow}(T^{(\la)}(x)) ]^{(\la)}$  is not covered by the
union $\Psi^{(\la)}( T^{(\la)}(\X \setminus x)), \la \in [\la_0, \infty).$  Similar germ-grain models have been considered in  section 4 of \cite{SY} and sections 2 and 4 of \cite{CSY}).
{\em We say that $T^{(\la)}(x)$ is extreme in $T^{(\la)}(\X)$ if $x \in \X$ is a vertex of $\rm{co}(\X)$.  Given $T^{(\la)}(\X) = \X'$, write $\rm{Ext}^{(\la)}(\X')$ for the set of extreme points in $\X'$. }


For $x_0 \in \R^d$ consider the half-space
$$
H(x_0):= \{x \in \R^d: \ \langle x, \frac{x_0} {|x_0|} \rangle \geq |x_0| \}.
$$
Taking again $\theta= d_{\S^{d-1}} ( \frac{x} {|x|}, \frac{x_0} {|x_0|} )$, we  rewrite $H(x_0)$ as
$$
H(x_0):= \{x \in \R^d: \ R_\la^2(1 - {|x_0| \over R_\la \cos \theta}) \geq R_\la^2(1 - {|x| \over R_\la}) \},
$$
Taking $T^{(\la)}(x_0)= (v_0, h_0)$ and using \eqref{defvh}, we see that $T^{(\la)}$ transforms $H(x_0)$ into the downward grain
\be \label{formu4.3}  T^{(\la)}(H(x_0)):=
[\Pi^{\downarrow}(v_0, h_0)]^{(\la)} := \{
(v,h) \in W_\la, \;  h \leq R_\la^2 - \frac{R_\la^2 - h_0} {
\cos[e_\la(v, v_0)]  }
 \}. \ee
Noting that $\R^d \setminus {\rm{co}}(\X)$ is the union of  half-spaces not containing
points in $\X$, it follows that $T^{(\la)}$ transforms $\R^d \setminus {\rm{co}}(\X)$ into the subset of $W_\la$
given by
$$
\Phi^{(\la) }(T^{(\la)}(\X)):=
\bigcup_{\left\{ \substack{w\in W_\la \\
[\Pi^{\downarrow}(w)]^{(\la)} \cap T^{(\la)}(\X)=\emptyset}\right.} [\Pi^{\downarrow}(w)]^{(\la)}.
$$
Thus $T^{(\la)}$  sends the boundary of ${\rm{co}}(\X)$ to the
continuous function on $W_\la$ 
whose graph coincides with the  boundary of $\Phi^{(\la)
}(T^{(\la)}(\X)).
$
There is thus a $1-1$ correspondence between convex hull boundaries
and a subset of the continuous functions on $\R^{d-1} \times \R$.
This contrasts with Eddy \cite{Eddy}, who mapped {\em support
functions} of convex hulls into a subset of the continuous functions
on $\R^{d-1} \times \R$.

The germ-grain models $\Psi^{(\la) }(\P^{(\la)} )$ and
$\Phi^{(\la) }(\P^{(\la)})$ link the geometry of
$K_\la$ with that of the limit paraboloid germ-grain models $\Psi(\P)$ and
  $\Phi(\P)$. Theorem \ref{Th4} and the upcoming Proposition \ref{BLem2} show that
  the boundaries  $\partial \Psi^{(\la) }(\P^{(\la)})$ and
$\partial \Phi^{(\la) }(\P^{(\la)})$ respectively converge in probability to  $\partial (\Psi(\P))$ and to  $\partial (\Phi(\P))$ as $\la \to \infty$.

The next lemma is suggestive of this convergence and shows for fixed
$w \in W_\la$ that $[\Pi^{\uparrow}(w)]^{(\la)}$ and
$[\Pi^{\downarrow}(w)]^{(\la)}$ locally approximate the paraboloids
$[\Pi^{\uparrow}(w)]^{(\infty)}:= \Pi^{\uparrow}(w)$ and
$[\Pi^{\downarrow}(w)]^{(\infty)}:= \Pi^{\downarrow}(w),$ respectively.  We may
henceforth refer to $[\Pi^{\uparrow}(w)]^{(\la)}$ and
$[\Pi^{\downarrow}(w)]^{(\la)}$ as {\em quasi-paraboloids} or sometimes `paraboloids' for short.
Recalling that $B_{d-1}(v,r)$ is the $(d-1)$ dimensional ball
 centered at $v \in \R^{d-1}$ with radius $r$, define the cylinder $C(v,r) \subset \R^{d-1} \times \R$ by
 \be \label{cyl}
 C(v,r):= C_{d-1}(v,r):=B_{d-1}(v,r) \times \R. \ee
{\em Here and in the sequel, by $c$ and $c_1, c_2,...$ we mean generic positive constants which may change from line to line. }

\begin{lemm} \label{paralem}  For all $w_1:= (v_1, h_1) \in W_\la$,  $L
\in (0, \infty)$, and all $\la \in [\la_0, \infty)$, we have \be \label{paralem-1}|| \partial(
[\Pi^{\uparrow}(w_1)]^{(\la)} \cap C(v_1, L) ) -
\partial( [\Pi^{\uparrow}(w_1)]^{(\infty)} \cap C(v_1, L)
)||_{\infty} \leq cL^3 R_\la^{-1} + ch_1 L^2 R_\la^{-2} \ee and \be
\label{paralem-2} || \partial( [\Pi^{\downarrow}(w_1)]^{(\la)} \cap
C(v_1, L) ) -
\partial( [\Pi^{\downarrow}(w_1)]^{(\infty)} \cap C(v_1, L)
)||_{\infty} \leq cL^3 R_\la^{-1} + ch_1 L^2 R_\la^{-2}. \ee
\end{lemm}

\noindent{\em Proof.} We first prove \eqref{paralem-1}.
By  \eqref{formu4.2} we have \be
\label{formu4.2-a}
\partial([\Pi^{\uparrow}(w_1)]^{(\la)}) := \{ (v,h) \in W_\la, \;
     h = R_\la^2 (1-\cos[e_\la(v, v_1)]) + h_1 \cos[e_\la(v,v_1)] \}.
     \ee
For $v \in B_{d-1}(v_1, L)$, notice that \be \label{ela} e_\la(v,
v_1) = |R_\la^{-1}v - R_\la^{-1}v_1| + O(|R_\la^{-1}v -
R_\la^{-1}v_1|^2) \ee and thus
$$
1 - \cos(e_\la(v, v_1)) = \frac{ |R_\la^{-1}v - R_\la^{-1}v_1|^2}
{2} + O(L^3 R_\la^{-3}).
$$
It follows that
$$
R_\la^2( 1 - \cos(e_\la(v, v_1))) = \frac{|v - v_1|^2} {2} + O(L^3
R_\la^{-1})
$$
and
$$
|h_1( 1 - \cos(e_\la(v, v_1)))| = O(h_1 L^2R_\la^{-2}).
$$
Thus the boundary of $[\Pi^{\uparrow}(w_1)]^{(\la)} \cap C(v_1, L)$ is
within $cL^3 R_\la^{-1} + ch_1 L^2 R_\la^{-2}$ of the graph of
$$
v \mapsto h_1 + \frac{ |v - v_1|^2}{2},
$$
which establishes \eqref{paralem-1}.  The proof of \eqref{paralem-2}
is similar, and goes as follows. By  \eqref{formu4.3} we have \be \label{formu4.3-a}
\partial([\Pi^{\downarrow}(w_1)]^{(\la)}) :=
     \{ (v,h) \in W_\la,
\;  h = R_\la^2 - \frac{R_\la^2 - h_1} { \cos[e_\la(v, v_1)]  }
 \}. \ee
Using \eqref{ela}, Taylor expanding $\cos\theta$ up to second order, and writing $1/(1 -r) = 1 + r + r^2 + ...$ gives
\be \label{formu4.3-b}
\partial( [\Pi^{\downarrow}(w_1)]^{(\la)} ) :=
     \{ (v,h) \in W_\la,
\;  h =  h_1 - \frac{ |v - v_1|^2} {2} +  O( R_\la^{-1} |v - v_1|^3 ) +
O( h_1 R_\la^{-2} |v - v_1|^2)    \}, \ee and \eqref{paralem-2}
follows.  \qed
\vskip.5cm

\noindent{\bf 3.2. The weak limit of $T^{(\la)}(\P_\la)$}. Put
$$\P^{(\la)}:= T^{(\la)}(\P_\la).$$
 Unlike the set-up of
 \cite{CSY}, the weak limit of $T^{(\la)}(\P_\la)$ converges to a point process which is non-homogenous and which is carried by all of $\R^{d-1} \times \R$.
 Let $\Vol_d$ denote
$d$-dimensional volume measure on $\R^d$ and let ${\rm{Vol}}_d^{(\la)}$ be the image of $R_\la\rm{Vol}_d$ under $T^{(\la)}$.
Recall the definition of $\P$ at \eqref{defP}.

\begin{lemm} \label{weakcon}  As $\la \to \infty$, we have
\  (a) $\P^{(\la)} \tod \P$  and  (b)
${\rm{Vol}}_d^{(\la)}
 \tod  \Vol_d.$
The convergence is in the sense of total variation convergence on
compact sets.
\end{lemm}


\noindent{\em Remarks.}  (i)  It is likewise the case that the image
of the binomial point process $\sum_{x \in \X_n} \delta_x$ under
$T^{(n)}$ converges in distribution to $\P$ as $n \to \infty$.

\noindent(ii) $T^{(\la)}$  carries $\P_\la$ into a
point process on $\R^{d-1} \times \R$ which in the large $\la$ limit is stationary
in the spatial coordinate.  This contrasts with the
transformation of Eddy \cite{Eddy} (and generalized in Eddy and Gale
\cite{EG}) which carries $\sum_{x \in \X_n} \delta_x$ into a point
process $(T_k, Z_k)$ on $\R \times \R^{d-1}$ where $T_k, k \geq 1,$ are points of
a Poisson point process on $\R$ with intensity $e^{-h} dh$ and
 $Z_k, k \geq 1,$ are i.i.d. standard Gaussian on $\R^{d-1}$.

\vskip.5cm \noindent{\em Proof}. Representing $x \in \R^d$ by $x =
ur, u \in \S^{d-1}, r \in [0, \infty)$,
 we find the image by $T^{(\la)}$ of the Poisson
measure on $\R^d$ with intensity
\be \label{Poisdens}
\la \phi(x) dx = \la \phi(ur) r^{d-1} dr d \sigma_{d-1}(u).
\ee
Make the change of variables
$$
v:= R_\la \exp_{d-1}^{-1}(u) = R_\la v_u, \ \  h := R_\la^2(1 - \frac{r} {R_\la}),
$$
The exponential map $\exp_{d-1} : T_{u_0} \S^{d-1}\to \S^{d-1}$ has the
following expression: 
\begin{equation}
  \label{expomap}
\exp_{d-1}(v)= \cos(|v|)(0,\cdots,0,1) +
\sin(|v|)(\frac{v}{|v|},0),\quad v\in\R^{d-1}\setminus\{\0\}.
\end{equation}
Therefore, since $v_u:= \exp_{d-1}^{-1}(u)$ we have
 $$d\sigma_{d-1}(u)= \sin^{d-2}(|v_u|) d(|v_u|) d\sigma_{d-2}(
 \frac{v_u} {|v_u|}) =   \frac{\sin^{d-2}(|v_u|) dv_u} {|v_u|^{d-2}}.$$
 Since $v_u = R_\la^{-1} v$, this gives
\be \label{3a}
 d\sigma_{d-1}(u)= \frac{ \sin^{d-2}( R_\la^{-1} |v| )}{ |R_\la^{-1} v|^{d-2} } (R_\la^{-1})^{d-1}dv.
 \ee We also have
\be \label{3b} r^{d-1} dr = [R_\la(1 - {h \over R_\la^2} )]^{d-1}
R_\la^{-1} dh \ee as well as
 \be \label{3c}
\la \phi(x) =  \la \phi(u  R_\la (1 - {h \over R_\la^2})) =
\sqrt{2  \log \la}
\exp\left(h - \frac{h^2}{2R_{\la}^2}\right). \ee Combining \eqref{Poisdens} and \eqref{3a}-\eqref{3c}, we get that $\P^{(\la)}$ has  intensity density \be
 \label{3d} { d \P^{(\la)}
\over dv dh} ((v,h)) = { \sqrt{2 \log \la} \over R_\la } \  \frac{ \sin^{d-2}( R_\la^{-1} |v| )}{ |R_\la^{-1} v|^{d-2} } (1 -
{h \over R_\la^2} )^{d-1} \exp(h - \frac{h^2}{2R_{\la}^2}), \ (v,h) \in W_\la.\ee
Given a fixed compact subset $D$ of $W_\la$,
this intensity
converges to the intensity of $\P$ in $L^1(D)$, completing the proof of part (a).

Replacing the intensity $\la \phi dx$ with $ dx$
in the above computations gives \be
 \label{Leb} { d \Vol_d^{(\la)}
\over dv dh} ((v,h)) = \frac{ \sin^{d-2}( R_\la^{-1} |v| )}{
|R_\la^{-1} v|^{d-2} } (1 - {h \over R_\la^2} )^{d-1}, \ (v,h) \in W_\la.\ee This intensity density
converges pointwise to $1$ as $\la \to \infty$, showing part (b).
 \qed

\vskip.5cm

\noindent{\bf 3.3. Re-scaled $k$-face and volume functionals}. Fix $\la \in [\la_0, \infty)$. {\em Let $\xi_k, k \in
\{0,1,...,d-1\},$ be a generic $k$-face functional,} as in
Definition \ref{kface}.
 The inverse transformation
$[T^{(\la)}]^{-1}$ defines generic re-scaled
functionals
$\txi$ defined for $\xi\in \Xi$, $w \in W_\la$ and $\X \subset W_\la$ by \be
\label{txidef} \txi(w, \X):= \xi^{(\la)}(w, \X):=
\xi([T^{(\la)}]^{-1}(w), [T^{(\la)}]^{-1}(\X)). \ee   For all 
$\la \in [\la_0, \infty),$ it follows that $
 \xi(x,
\P_\la):= \xi^{(\la)} (T^{(\la)}(x), \P^{(\la)}).$ Note for all
$\la \in [\la_0, \infty)$, $k \in \{0,1,...,d-1\},$  $w_1 \in W_\la$, and $\X \subset W_\la$, that $\xi^{(\la)}_k(w_1, \X)$
is the product of $(k + 1)^{-1}$ and the number of quasi-parabolic
$k$-dimensional faces of $\partial (\bigcup_{w \in \X}
[\Pi^{\downarrow}(w)]^{(\la)})$ which contain $w_1$, $w_1 \in \rm{Ext}^{(\la)}(\X)$, otherwise $\xi^{(\la)}_k(w_1, \X)= 0$.

\vskip.2cm

Similarly, define for  $w \in \rm{Ext}^{(\la)}(\X)$
\be \label{re-vol}
 \xi^{(\la)}_V(w, \X) =  {1 \over d} \int_{ v\in \rm{Cyl}^{(\lambda)}(w) } \int_0^{ \Phi^{(\la)}(\X)(v)}
 d {\rm{Vol}}_d^{(\la)}((v,h))), \ee
where ${\rm{Cyl}}^{(\lambda)}(w):= {\rm{Cyl}}^{(\lambda)}(w, \X)$ denotes the projection onto $\R^{d-1}$ of the quasi-parabolic faces of $\Phi^{(\la)}(\X)$ containing $w$.  When $w \notin \rm{Ext}^{(\la)}(\X)$ we
define $\xi^{(\la)}_V(w, \X) = 0$.

Given $\la \in [\la_0, \infty),$ let $\Xi^{(\la)}$ denote the collection of re-scaled functionals
$\xi^{(\la)}_{k}, k \in \{0,1,...,d-1\},$ together with
$\xi^{(\la)}_{V}$.
Our main goal in the next section is to show
that, given a generic  $\xi^{(\la)} \in
\Xi^{(\la)}$, the means and covariances of $\xi^{(\la)}(\cdot,
\P^{(\la)})$ converge as $\la \to \infty$ to the respective means
and covariances of $\xi^{(\infty)}(\cdot, \P)$, with
$\xi^{(\infty)} \in \Xi^{(\infty)}$.

\vskip.1cm






\vskip.3cm

\section{Properties of  re-scaled $k$-face and volume functionals} \label{Sec4}

\allco

To establish convergence of re-scaled functionals $\xi^{(\la)} \in \Xi^{(\la)},
\la \in [\la_0, \infty),$  to their respective counterparts $\xi^{(\infty)} \in
\Xi^{(\infty)}$, we first need to show that $\xi^{(\la)} \in
\Xi^{(\la)}, \la \in [\la_0, \infty]$  satisfy a localization in the spatial and time coordinates
$v$ and $h$, respectively.  These
localization results are  the analogs of Lemmas 7.2 and 7.3 of \cite{CSY}.
{\em In the following the point process $\P^{(\la)}$, $\la = \infty$, is taken to be
$\P$ whereas  $W_\la$, $\la = \infty$, is taken to be $\R^d$.} Many of our proofs for the case $\la \in (0, \infty)$ may be modified to yield explicit proofs of some
unproved assertions in \cite{CSY}.

\vskip.5cm \noindent{\bf 4.1. Localization of $\txi$.} 
Recall the definition at \eqref{cyl} of the cylinder $C(v,r):= C_{d-1}(v,r):=B_{d-1}(v,r) \times \R$, $v\in \R^{d-1}$, $r>0$.
Given a generic functional $\xi^{(\la)} \in \Xi^{(\la)}, \la \in [\la_0, \infty],$ and $w:=(v,h) \in W_\la$, we shall
 write \be \label{cyl-1} \xi^{(\la)}_{[r]}(w,\P^{(\la)}) := \xi^{(\la)}(w, \P^{(\la)} \cap  C_{d-1}(v,r)).\ee

Given $\xi^{(\la)}, \la \in [\la_0, \infty],$ recall from \cite{CSY,SY} that a random variable $R :=
R^{\xi^{(\la)}}[w]:=R^{\xi^{(\la)}}[w, \P^{(\la)}]$ is a
 {\it spatial localization radius} for  $\xi^{(\la)}$ at $w$ with respect to $\P^{(\la)}$ iff a.s.
 \be \label{slr} \xi^{(\la)}(w,\P^{(\la)}) = \xi^{(\la)}_{[r]}(w,\P^{(\la)}) \mbox{ for all } r \geq R. \ee
 There are in general more than one $R$ satisfying \eqref{slr} and we shall henceforth assume $R$ is the infimum of all reals
satisfying \eqref{slr}.  

We may similarly define  a localization radius in the non-rescaled picture. Indeed, given a generic functional $\xi$ and $x\in\R^d\setminus \{0\}$, we shall
 write
 \begin{equation*}
\xi_{[r]}(x,\P_{\la}) := \xi(x, \P_{\la} \cap S(x,r) ).
 \end{equation*}
where $S(x,r):=\{y\in\R^d\setminus\{0\}: d_{{\mathbb S}^{d-1}}(x/|x|,y/|y|)\le r\}$. $R^{\xi}[x,\P_\la]$ is then the infimum of all $R \in (0, \infty)$ which satisfy $\xi(x,\P_{\la}) = \xi_{[r]}(x,\P_{\la})$ for every $r \in [R, \infty)$.
In particular, by rotation-invariance of $\P_{\la}$ and the fact that $|v-\0|=d_{{\mathbb S}^{d-1}}(\exp_{d-1}(\0),\exp_{d-1}(v))$ for all $v\in \B_{d-1}(\pi)$, we have the distributional equality:
\begin{equation}\label{locradequal}
R^{\xi}[x,\P_\la]\overset{D}{=}R^{\xi}[|x|u_0,\P_\la]
=R_{\la}^{-1}R^{\xi^{(\la)}}[(\0,h_0)],
\end{equation}
where $h_0=R_\la^2\left(1-{|x|}/{R_\la}\right)$.
In view of \eqref{locradequal}, it is enough to investigate the distribution tail of $R^{\xi^{(\la)}}[(\0,h_0)]$ for any $h_0\in\R$. In the next lemmas, we prove that
the functionals
$\xi^{(\la)} \in \Xi^{(\la)}, \la \in [\la_0, \infty],$
  admit spatial localization radii with tails decaying
 super-exponentially fast at $(\0,h_0)$, $h_0\in\R$.
We first establish a localization radius for $\xi_0$.  We remark this shows that
$\rm{Ext}^{(\la)}(\P^{(\la)}),  \la \in [\la_0, \infty],$ is a strongly mixing
random point set.

\begin{lemm} \label{lem4.0}
There is a constant
 $c>0$ such that the localization radius  $R^{\xi_0^{(\la)}}[(\0,h_0)]$
 satisfies for all $\la \in [\la_0, \infty]$, $h_0 \in (-\infty, R_\la^2]$, and $t\ge (-h_0\vee 0) $
  \begin{equation}\label{LocalExpr1-0}
   P[R^{\xi_0^{(\la)}}[(\0,h_0)] >  t] \leq c \exp(- \frac{t^{2} }{ c}).
\end{equation}
\end{lemm}

\noindent{\em Proof.}
Abbreviate $\xi_0$ by $\xi$.
It suffices to show that \eqref{LocalExpr1-0} holds for
$t \geq -h_0 \vee c$, $c$ a positive constant, a simplification used repeatedly in what follows. 
For $t\ge (-h_0\vee 0)$ and $\la\in [\la_0,\infty]$, we have
\begin{equation}
  \label{eq:decompE1E2}
\{R^{\xi^{(\la)}}[(\0,h_0)] >  t\}\subset E_1\cup E_2,
\end{equation}
where
$$ E_1:=  \left\{ \mbox{$R^{\xi^{(\la)}} [(\0,h_0)] >  t, \ (\0,h_0) \notin {\rm{Ext}^{(\la)}}(\P^{(\la)})$} \right\}$$ and $$E_2:= \left\{ \mbox{$R^{\xi^{(\la)}} [(\0,h_0)] >  t, \ (\0,h_0) \in {\rm{Ext}^{(\la)}}(\P^{(\la)})$} \right\}.$$
Rewrite $E_1$ as
$$
E_1= \{(\0,h_0) \notin {\rm{Ext}^{(\la)}}(\P^{(\la)}), \ (\0,h_0) \in {\rm{Ext}^{(\la)}}(\P^{(\la)}\cap C(\0,t)) \}.
$$
If $E_1$ occurs then there is a
$$
w_1:=(v_1,h_1) \in \partial\left( [\Pi^{\uparrow}((\0,h_0))]^{(\la)} \right) \cap C(\0, t)
$$
belonging to some $[\Pi^{\uparrow}(y)]^{(\la)}, \ y \in \P^{(\la)} \cap C(\0,
t)^c$, but $w_1 \notin \bigcup_{w \in  \P^{(\la)} \cap C(\0, t) }
[\Pi^{\uparrow}(w)]^{(\la)}.$ In other words,  $w_1$ is covered by paraboloids
with apices in $\P^{(\la)}$, but not by paraboloids with apices in $\P^{(\la)} \cap C(\0, t)$.  This means that the down paraboloid
$[\Pi^{\downarrow}(w_1)]^{(\la)}$
does not contain points in  $C(\0, t) \cap \P^{(\la)}$, but it
 must contain a point in $C(\0, t)^c \cap \P^{(\la)}$.
In other words, we have $E_1 \subset F_1 \cup F_2$, where
\begin{align*}
F_1:= & \{ \exists w_1:=(v_1,h_1) \in \partial [\Pi^{\uparrow}((\0,h_0))]^{(\la)} \cap C(\0, t): \ h_1 \in (- \infty, t),
\\ & \ \ [\Pi^{\downarrow}(w_1)]^{(\la)} \cap
C(\0, t) \cap \P^{(\la)} = \emptyset \ , [\Pi^{\downarrow}(w_1)]^{(\la)} \cap
C(\0, t)^c \cap \P^{(\la)} \neq \emptyset \}
\end{align*}
and
\begin{align*}
F_2:= & \{ \exists w_1:=(v_1,h_1) \in \partial[\Pi^{\uparrow}((\0,h_0))]^{(\la)} \cap C(\0, t):  \ h_1 \in [t,  \infty),\\
& \ \  [\Pi^{\downarrow}(w_1)]^{(\la)} \cap
C(\0, t) \cap \P^{(\la)} = \emptyset \}.
\end{align*}

If $E_2$ happens then there is $w_1:= (v_1, h_1) \in C(\0,t)^c \cap [\Pi^{\uparrow}((\0,h_0))]^{(\la)}$
which is not covered by paraboloids with apices in $\P^{(\la)}$ and $(\0,h_0)$ belongs to $\partial [\Pi^{\downarrow}(w_1)]^{(\la)}$.
  Notice that $w_1\in C(\0,\pi R_\la/2)$ since the ball $[T^{(\la)}]^{-1}([\Pi^{\uparrow}((\0,h_0))]^{(\la)})$ is included in $\R^{d-1} \times [0, \infty)$.  There is a constant $c>0$ such that $1-\cos(\theta)\ge c \theta^2$ for $\theta\in [0,\pi]$ so that in view of \eqref{formu4.2} and $e_\la(v_1, \0) = R_\la^{-1} |v_1|$, we have
$$h_1\ge h_0\cos(R_{\la}^{-1}|v_1|)+R_\la^2(1-\cos(R_{\la}^{-1}|v_1|))\ge (h_0\wedge 0) + c|v_1|^2 \ge (h_0\wedge 0) + ct^2.$$
Now $h_0\wedge 0 \geq -t$ always holds so
we obtain $h_1 \geq -t + ct^2 \geq t$ for large enough $t$.
Thus we have $E_2 \subset \tilde{E}_2 $ where
\be \label{E2}
\tilde{E}_2:=\{ \exists w_1:=(v_1,h_1) \in \partial[\Pi^{\uparrow}((\0,h_0))]^{(\la)}: \ h_1 \in [t, \infty),
[\Pi^{\downarrow}(w_1)]^{(\la)} \cap \P^{(\la)} = \emptyset \}.
\ee
By  \eqref{eq:decompE1E2} and the inclusions $E_1\subset F_1\cup F_2$ and $E_2\subset \tilde{E}_2$,
it is enough to show that each term $P[F_1]$, $P[F_2]$ and $P[\tilde{E}_2]$ is bounded by $c \exp(-t^2/c)$.




\vskip.5cm

\noindent{\it Upper-bound for $P[F_1]$}.
We start with the case $\lambda=\infty$.
Consider a fixed $w_1\in \partial \Pi^{\uparrow}((\0,h_0))$ with
$h_1=h_0+\frac{1}{2}|v_1|^2\le t$. The probability that
$\Pi^{\downarrow}(w_1)\cap C(\0,t)^c\cap \P\ne \emptyset$ is bounded by
the  $d\P$ measure of $\Pi^{\downarrow}(w_1)\cap C(\0,t)^c$. The maximal height of  $\Pi^{\downarrow}(w_1)\cap C(\0,t)^c$ is $h_1-\frac{1}{2}(t-\sqrt{2(h_1-h_0)})^2$.
Consequently, the $d\P$ measure of  $\Pi^{\downarrow}(w_1)\cap
C(\0,t)^c$ is bounded by the $d\P$ measure of
$\Pi^{\downarrow}(w_1)\cap \{(v,h):h\le
h_1-\frac{1}{2}(t-\sqrt{2(h_1-h_0)})^2\}$. Recall that $c$ is a
constant which changes from line to line. Up to a multiplicative constant, the  $d\P$ measure of $\Pi^{\downarrow}(w_1)\cap C(\0,t)^c$ is bounded by
\begin{align*}
\int_{-\infty}^{h_1-\frac{1}{2}(t-\sqrt{2(h_1-h_0)})^2}
e^{h}(2(h_1-h))^{\frac{d-1}{2}}dh&=
e^{h_1}\int_{\frac{1}{2}(t-\sqrt{2(h_1-h_0)})^2}^{+\infty} e^{-u}(2u)^{\frac{d-1}{2}}du\\
&\le c \exp\left(h_1- {1 \over c }(\frac{t^2}{2}+h_1-h_0-t\sqrt{2(h_1-h_0)})\right),
\end{align*}
where we put $u := h_1 - h$.

Consequently, discretizing
$\partial \Pi^{\uparrow}((\0,h_0))\cap (\R^{d-1}\times (h_0, t])$ and using $h_0\le h_1\le t$,  we get
\begin{align*}
P[F_1] &\leq ce^{ {-t^2 \over 2c }}\int_{h_0}^{t}(2(h_1-h_0))^{\frac{d-2}{2}}\exp \left((1-{1\over c})h_1+ {h_0 \over c}+ {t\sqrt{h_1-h_0} \over c} \right) dh_1\\
& \le ce^{ {-t^2 \over 2c }} \int_{0}^{t-h_0}(2h_1)^{\frac{d-2}{2}}\exp \left((1-{1\over c})h_1
+h_0
+ {t\sqrt{h_1} \over c} \right)dh_1\\
& \le c e^{ { -(t^2-t^{3/2}-t) \over  c}}\\
& \le c e^{ {-t^2 \over c}},
\end{align*}
concluding the case $\la = \infty$.

When $\la \in [\lambda_0,\infty)$,
%
%
 recall from \eqref{formu4.2} that
\be \label{4.2aa}
[\Pi^{\uparrow}((\0,h_0))]^{(\la)} := \{ (v,h) \in W_\la, \;
h \geq  R_\la^2 (1-\cos[e_\la(v, \0)]) + h_0 \cos[e_\la(v,\0)] \}.
\ee
We claim that
$[\Pi^{\uparrow}((\0,h_0))]^{(\la)} \cap (\R^{d-1} \times (-\infty, t])$ has a spatial diameter (in the $v$ coordinates)  bounded by $c_1 \sqrt{t}$.   We see this as follows. Let $(v, h) \in [\Pi^{\uparrow}((\0,h_0))]^{(\la)} \cap (\R^{d-1}\times(-\infty, t])$.  When $h \leq t$ and $|h_0| \leq t$, the display \eqref{4.2aa} yields $R_\la^2 (1-\cos[e_\la(v, \0)]) \leq 2t$.   Thus $1-\cos[e_\la(v, \0)] \leq 2t R_\la^{-2}$. It follows that
\begin{equation}\label{eq:apertquasiparab}
c e_\la(v, \0)^2  \leq   1-\cos[e_\la(v,\0)] \leq 2t R_\la^{-2}.
\end{equation}
Using the equality
$e_\la(v, \0) = R_\la^{-1}|v|$, we deduce $|v| \leq c_1 \sqrt{t}$, as desired.

Let
$$
w_1:=(v_1,h_1) \in \partial [\Pi^{\uparrow}((\0,h_0)) ]^{(\la)} \cap
C(\0, t).
$$
We now estimate the maximal height of $[\Pi^{\downarrow}(w_1) ]^{(\la)} \cap C(\0, t)^c.$
If $(v,h)$ belongs to the boundary of $[ \Pi^{\downarrow}(w_1) ]^{(\la)}$ then we have from \eqref{formu4.3-a} that
$$
h = R_\la^2 - \frac{R_\la^2 - h_1} { \cos[e_\la(v, v_1)]  }
$$
which gives
\begin{equation}
  \label{eq:bounddist}
c e_\la(v, v_1)^2 \leq 1-\cos[e_\la(v, v_1)] =\frac{h_1 - h} {R_\la^2 - h} \leq \frac{h_1 - h}  {R_\la^2 - t} \leq \frac{h_1 - h} {R_\la^2 -
2\pi R_\la}
\end{equation}
where we use $h \leq t \leq 2\pi R_\la$. Indeed, we may without loss of generality assume  $t \in [0, 2\pi R_\la]$, since the stabilization radius never exceeds the spatial diameter of $W_\la$.
Consequently, we have
\be \label{E2}
h \leq h_1 - c_2(R_{\la}^2-2\pi R_{\la})e_\la(v, v_1)^2.
\ee
The maximal height of  $[( \Pi^{\downarrow}(w_1) ]^{(\la)} \cap C(\0, t)^c$ is found by letting $v$ belong to  the boundary of $C(\0,t)$. In particular, we have $e_{\la}(v,\0)= R_{\la}^{-1}|v|= R_{\la}^{-1}t$. Moreover, we deduce from \eqref{eq:apertquasiparab} that $e_{\la}(v_1,\0)\le c_1 R_{\la}^{-1}\sqrt{t}$. Consequently, we have
$$e_{\la}(v,v_1)\ge e_{\la}(v,\0)-e_{\la}(v_1,\0)\ge R_{\la}^{-1}(t-c_1\sqrt{t}).$$
Combining the last inequality above with \eqref{E2} shows  for any
$(v,h) \in \partial[\Pi^{\downarrow}(w_1) ]^{(\la)}\cap \partial C(\0,t)$  that
$$h\le h_1-c_2\frac{R_{\la}^2-2\pi R_{\la}}{R_{\la}^2}(t-c_1\sqrt{t})^2\le t - c_3 (t - c_1 \sqrt{t})^2.$$

Now we follow the proof for the case $\la = \infty$.  We have
$$d\P^{(\la)}( [\Pi^{\downarrow}(w_1) ]^{(\la)} \cap C(\0, t)^c ) \leq d\P^{(\la)}
([\Pi^{\downarrow}(w_1) ]^{(\la)} \cap \{(v,h): \ h \leq t - c_2 (t - c_1 \sqrt{t})^2 \}),
$$
In view of \eqref{3d} and \eqref{eq:bounddist}, $d\P^{(\la)}( [\Pi^{\downarrow}(w_1) ]^{(\la)} \cap C(\0, t)^c )$ is bounded by
$$
c \int_{-\infty}^{ t - c_3 (t - c_1 \sqrt{t})^2 }
(1 - {h \over R_\la^2} )^{d-1} e^{h}
\left[\int {\bf 1}(e_\la(v,v_1)\le {c\sqrt{h_1-h} \over R_\la})\frac{ \sin^{d-2}( R_\la^{-1} |v| )}{ |R_\la^{-1} v|^{d-2} }dv\right]dh.
$$
Using the change of variables $u=\exp_{d-1}^{-1}(R_\la^{-1}v)$ with $u_1=\exp_{d-1}^{-1}(R_\la^{-1}v_1)$ gives
 \begin{align}\label{covsphere}
& d\P^{(\la)}([\Pi^{\downarrow}(w_1) ]^{(\la)} \cap C(\0, t)^c) \nonumber \\&\le   c \int_{-\infty}^{ t - c_3 (t - c_1 \sqrt{t})^2 } (1 - {h \over R_\la^2} )^{d-1} e^{h}
\left[\int_{{\mathbb S}^{d-1}}
{\bf 1}(d_{{\mathbb S}^{d-1}}(u,u_1)\le { c\sqrt{h_1-h} \over R_\la} )
R_\la^{d-1}d\sigma_{d-1}(u)\right]dh
\nonumber \\
&\le  c \int_{-\infty}^{ t - c_3 (t - c_1 \sqrt{t})^2 } (1 - {h \over R_\la^2} )^{d-1} e^{h} (h_1-h)^{(d-1)/2}dh.
 \end{align}
For $t$ large the upper limit of integration is at most $-c_4 t^2$ where $c_4=c_3/2$. 
 There is a positive constant $c_5$ such that $(1 - {h/ R_\la^2} )^{d-1} e^{h/2} \leq e^{c_5 h}$ holds for all $h \in (-\infty, 0]$.
Also, $$(h_1 - h)^{(d-1)/2 } \leq c(t^{(d-1)/2 } + |h|^{(d-1)/2 }). $$  Putting these estimates together yields
$$d\P^{(\la)} ([\Pi^{\downarrow}(w_1) ]^{(\la)} \cap C(\0, t)^c) \leq
c \int_{-\infty}^{ -c_4 t^2 }  e^{c_5 h} (t^{(d-1)/2 } + |h|^{(d-1)/2 }) dh \leq c_6 \exp(- t^2/c_6).
$$
%
Consequently, discretizing
$\partial [\Pi^{\uparrow}((\0,h_0))]^{(\la)}  \cap C(\0,t)\cap(\R^{d-1}\times (-\infty, t])$
we get for $\la \in [\la_0, \infty)$
\begin{align*}
P[F_1] &\leq c_6 e^{-t^2/c_6}\int_{h_0}^{t} t^{d-2}  dh_1 \le c_7 e^{-t^2/c_7}.\\
\end{align*}

\noindent{\em Upper-bound for $P[F_2]$}. 
 We again  start with the case $\la=\infty$. Suppose $h_1 \in [t, \infty)$ with $t$ large. As noted, $\Pi^{\downarrow} \cap  C(\0, t)$ does not contain
points in $\P$.
The $d\P$ measure of $\Pi^{\downarrow}(w_1) \cap C(\0, t)$
is bounded below by the $d\P$ measure of $\Pi^{\downarrow}(w_1) \cap C(\0, t)  \cap  (\R^{d-1} \times [0, \infty))$, which we
generously bound below by $e^{h_1/2}$.
Thus the probability that $\Pi^{\downarrow} \cap  C(\0, t)$ does not contain
points in $\P \cap C(\0, t)$ is bounded above by $\exp(- e^{h_1/2})$.

Discretizing $\partial \left( \Pi^{\uparrow}((\0,h_0)) \right)\cap (\R^{d-1} \times [t, \infty)) \cap C(\0, t)$ into unit cubes, we see that the probability that there is
$w_1:=(v_1,h_1) \in \partial \left( \Pi^{\uparrow}((\0,h_0)) \right) \cap C(\0, t)$
such that  $\Pi^{\downarrow}(w_1)$ does not contain
points in $\P \cap C(\0, t)$ is bounded by
$$
c \int_t^{\infty} t^{d -2} \exp(- e^{h_1/2}) dh_1 \leq c t^{d -2} \exp({-e^{t/2} \over c}).
$$
Thus there is a  constant $c$ such that
$P[F_2] \leq c \exp(- t^{2}/ c)$ for $t \geq (- h_0 \vee c)$.
~\\~\\
When $\la\in [\la_0,\infty)$, we proceed as follows.
Let $w_1$ be the point defined in event $F_2$. Let $S$ be the unit volume cube centered at $(v_1-\sqrt{d-1}v_1/(2|v_1|), (h_1 + 1)/2)$.
We claim that for $t$ large enough, $S$ is included in $[\Pi^{\downarrow}(w_1)]^{(\la)}\cap C(\0,t\wedge 3\pi R_\la/4)$. 
Indeed, $S$ is clearly included in $[\Pi^{\downarrow}(w_1)]^{(\la)}\cap C(\0,t)\cap (\R^{d-1}\times [h_1/4,\infty))$ and since $v_1\in B_{d-1}(\0, \pi R_\la/2)$, $S$ is included in $C(\0,3\pi R_\la/4)$.
By \eqref{3d} there exists  a constant $c>0$ such that for all
 $(v,h) \in S$
$$\frac{d\P^{(\la)}((v,h))}{dv dh}\ge c(1 -   \frac{ h_1}{4R_\la^2} )^{d-1}  \exp( \frac{h_1}{4}  - \frac{h_1^2}{32R_\la^2}).$$
Now using $h_1 \in (-\infty,  R_{\la}^2]$, we obtain  ${d\P^{(\la)}(S)}\ge c \exp (7h_1/32).$ Consequently, the probability that $ [\Pi^{\downarrow}(w_1)]^{(\la)}\cap C(\0,t)\cap \P^{(\la)}=\emptyset$ is bounded above by $\exp(-ce^{h_1/c})$. Discretizing
$\partial [\Pi^{\uparrow}((\0,h_0))]^{(\la)}  \cap (\R^{d-1} \times [t, \infty)) \cap C(\0, t)$, we see that the probability that there is
$w_1:=(v_1,h_1) \in \partial [\Pi^{\uparrow}((\0,h_0))]^{(\la)} \cap C(\0, t)$
such that  $[\Pi^{\downarrow}(w_1)]^{(\la)}$ does not contain
points in $\P^{(\la)} \cap C(\0, t)$ is bounded by
$$
c \int_t^{\infty} t^{d -2} \exp(- ce^{h_1/c}) dh_1 \leq c \exp(-ce^{t/c}).
$$
~\\~\\
\noindent{\em Upper-bound for $P[\tilde{E}_2]$}. The arguments closely follow those for $P[F_2]$ and we sketch the proof only for finite $\la$ as the case $\la=\infty$ is similar.
As above, consideration of the cube $S$  shows that
 $P[ [\Pi^{\downarrow}(w_1)]^{(\la)}\cap \P^{(\la)}=\emptyset]$ is bounded above by $\exp(-ce^{h_1/c})$. Only the discretization differs from the case of $P[F_2]$. Indeed, we need now to discretize $\partial [\Pi^{\uparrow}((\0,h_0))]^{(\la)}  \cap (\R^{d-1} \times [t, \infty))$. We use the fact that
  $|v_1| \le c R_\la \sqrt{h_1-h_0} / \sqrt{R_\la^2-h_0}$ as soon as $(v_1,h_1)\in \partial[\Pi^{\uparrow}((\0,h_0))]^{(\la)}$ (by \eqref{formu4.2} and the  arguments as in \eqref{eq:apertquasiparab}). We obtain
$$P[\tilde{E}_2]\le
c \int_{t\vee h_0}^{R_\la^2} [\frac{R_\la}{\sqrt{R_\la^2-h_0}}\sqrt{h_1-h_0}]^{d -2} \exp(- ce^{h_1/c}) dh_1.
$$
When $h_0 \in (-\infty, R_\la^2/2]$, we bound the ratio $R_\la /\sqrt{R_\la^2-h_0}$ by $\sqrt{2}$ to obtain
 $P[\tilde{E}_2]\le c \exp(-ce^{t/c})$. When $h_0 \in (R_\la^2/2, R_\la^2]$, we bound $(h_1-h_0)/ (R_\la^2-h_0)$ by $1$ and we bound $\exp(- ce^{h_1/c})$ by $\exp(- ce^{h_1/c}/2-ce^{R_\la^2/(2c)}/2)$ and we also obtain $P[\tilde{E}_2]\le c \exp(-ce^{t/c})$, as desired.
\vskip.3cm
Combining the above bounds for $P[F_1], P[F_2],$ and $P[\tilde{E}_2]$ thus yields
$$P[E_1] + P[E_2] \leq P[F_1] + P[F_2] +P[\tilde{E}_2]\leq c \exp(- {t^{2} \over c}),$$ showing Lemma \ref{lem4.0}
as desired. \qed

\vskip.5cm Whereas Lemma \ref{lem4.0} localizes $k$-face and volume functionals in
 the spatial domain, we now  localize in the height/time domain.  We show that  the boundaries of the paraboloid germ-grain processes $\Psi^{(\la)}(\P^{(\la)})$
and $\Phi^{(\la)}(\P^{(\la)})$, $\la \in [\la_0, \infty]$,  are not far from $\R^{d-1}$. Recall that $\P^{(\la)}, \la = \infty$, is taken to be $\P$ and we also write $\Psi(\P)$ for $\Psi^{(\infty)}(\P^{\infty})$.
If $w \in \rm{Ext}^{(\la)}(\P^{(\la)})$ we put $H(w):=H(w,\P^{(\la)})$ to be the maximal height coordinate (with respect to $\R^{d-1}$) of an apex of a down paraboloid which contains a parabolic face in $\Phi^{(\la)}(\P^{(\la)})$
containing $w$, otherwise we put $H(w) = 0$.

\begin{lemm}\label{BLem1}  (a) There is a constant $c$ such that for all $\la \in [\la_0, \infty]$, $h_0 \in (-\infty, R_{\la}^2],$ and $t  \in [h_0 \vee 0, \infty)$
we have
\be \label{fbd-1} P[H((\0,h_0), \P^{(\la)}) \geq t] \leq c \exp(-\frac{e^t} {c}). \ee
\noindent(b) There is a constant $c$ such that for all $L \in (0, \infty)$,  $t \in (0, \infty)$,  and  $\la \in [\la_0, \infty]$ 
we have
\be \label{sbd}
 P[ || \partial \Psi^{(\la)}(\P^{(\la)}) \cap
C(\0,L) ||_{\infty} > t] \leq cL^{2(d-1)} e^{- \frac{t} {c} }.
\ee
The bound \eqref{sbd} also holds for the dual process
$\Phi^{(\la)}(\P^{(\la)}).$
\end{lemm}

\noindent{\em Proof.} Let us first prove \eqref{fbd-1}. We do this for $\la \in [\la_0, \infty)$ and we claim that a  similar proof holds for $\lambda=\infty$. Rewrite the event $\{H((\0,h_0), \P^{(\la)}) \geq t\}$ as
\begin{align*}
\{H((\0,h_0), \P^{(\la)}) \geq t\}&=\{ \exists w_1:=(v_1,h_1) \in \partial[\Pi^{\uparrow}((\0,h_0))]^{(\la)}:  \ h_1 \in [t,  \infty),\\
& \ \ \hspace*{4.5cm}  [\Pi^{\downarrow}(w_1)]^{(\la)}  \cap \P^{(\la)} = \emptyset \}.
\end{align*}

Let us consider $w_1 :=(v_1,h_1) \in \partial[\Pi^{\uparrow}((\0,h_0))]^{(\la)}$ and
put  $[T^{(\la)}]^{-1}(\0,h_0) := \rho u_0$, $\rho\in [0, \infty)$.   Since $[\Pi^{\uparrow}((\0,h_0))]^{(\la)}$ is the image by $T^{(\la)}$ of the ball $B_{d}(\rho u_0/2,\rho/2)$, it is a subset of the image of the upper-half space, i.e. a subset of $C(\0,\pi R_\la/2)$. Consequently, the unit-volume cube centered at $(v_1,h_1-1)$ is included in $ [\Pi^{\downarrow}(w_1)]^{(\la)} \cap C(\0,3\pi R_\la/4)$. The proof now follows along the same lines as for the bound for $P[\tilde{E}_2]$ in the proof of Lemma \ref{lem4.0}. The ${\P}^{(\la)}$-measure of that cube exceeds $c\exp(h_1/c)$. The probability that $ [\Pi^{\downarrow}(w_1)]^{(\la)}  \cap \P^{(\la)} = \emptyset$ is bounded above by $c \exp(-ce^{h_1/c})$. Discretizing $(\R^{d-1} \times [t, \infty))$ into unit cubes, we obtain \eqref{fbd-1}.

We now prove \eqref{sbd}. We bound the probability of the two events
  $$E_3:=\{\partial\Psi^{(\la)}(\P^{(\la)})   \cap\{(v,h):|v|\le L,h>t\}\ne \emptyset\}$$ and $$E_4:=\{\partial\Psi^{(\la)}(\P^{(\la)}) \cap\{(v,h):|v|\le L,h<-t\}\ne \emptyset\}.$$

When in $E_3$, there is a point $w_1:= (v_1,h_1)$ with $ h_1 \in [t, \infty), |v_1| \leq L,$ and such that $[\Pi^{\downarrow}(w_1)]^{(\la)}\cap \P^{(\la)}=\emptyset$.
Following the proof of  \eqref{fbd-1}, we construct a unit-volume cube in $C(\0,L)$ which is a domain where the density of the $d\P^{(\la)}$  measure exceeds $ce^{h_1/c}$.  Discretization of $\{(v,h): \ |v| \leq L, h \in [t, \infty)\}$
into unit volume sub-cubes gives  $$P[E_3] \leq cL^{d-1} \exp(-{ e^t\over c}).$$

On the event $E_4$, there exists a point $(v_1,h_1)$ with $|v_1|\le L$ and $h_1 \in (-\infty, -t]$ which is on the boundary of an upward
paraboloid with apex in $\P^{(\la)}$. The apex of this upward paraboloid is
contained in the union of all down paraboloids with apex on
$B_{d-1}(\0,L)\times \{h_1\}$. The $d\P^{(\la)}$  measure of this union is bounded by $cL^{d-1} \exp(h_1/c).$
Consequently, the probability that the union contains
points from $\P^{(\la)}$ is less than $1-\exp(-cL^{d-1} e^{h_1/c}) \le cL^{d-1} \exp(h_1/c).$
It remains to discretize and integrate over
$h_1\in (-\infty,-t)$.  This goes as follows.

Discretizing $C(\0,L) \times (-\infty,-t]$ into unit volume subcubes
and using  \eqref{3d}, we find that the probability  there exists
$(v_1,h_1) \in \R^{d-1} \times (-\infty, -t]$ on the boundary of an up paraboloid
is thus bounded by
 $$
 cL^{2(d-1)}\int_{-\infty}^{-t}e^{h_1/c}(1-\frac{h}{R_{\la}^2})^{d-1}e^{h_1}dh_1
$$
 This establishes \eqref{sbd}.  The same argument applies to the dual process $\Phi^{(\la)}(\P^{(\la)})$. 
 \qed

\vskip.5cm

We now extend  Lemma \ref{lem4.0} to all $\xi \in \Xi$. 
\begin{lemm} \label{lem4.01}
There is a constant $c >0$ such that for all $\xi \in \Xi$, $\la \in [\la_0, \infty]$, and $h_0 \in (-\infty, R_\la^2]$, the localization radius  $R^{\xi^{(\la)}}[(\0,h_0)]$
 satisfies for all $t \in [|h_0|, \infty)$
  \begin{equation}\label{LocalExpr1}
   P[R^{\xi^{(\la)}}[(\0,h_0)] >  t] \leq c \exp(- \frac{t^{2} }{ c}).
\end{equation}
\end{lemm}
\noindent{\em Proof.} We show \eqref{LocalExpr1} for $\lambda \in [\la_0,\infty)$, as the proof is analogous for $\lambda= \infty$.
When $H((\0,h_0),\P^{(\la)})\le t$, then 
$\xi^{(\la)}((\0,h_0))$ only depends on points of $\P^{(\la)}$
in
$${\mathcal U}:=\bigcup_{w_1\in [\Pi^{\uparrow}((\0,h_0))]^{(\la)}\cap\R^{d-1}\times (-\infty,t]}[\Pi^{\downarrow}(w_1)]^{(\la)}.$$
Let $w=(v,h)\in {\mathcal U}$  and $w_1=(v_1,h_1)$, $h_1\le t$,  be such that $\partial [\Pi^{\downarrow}(w_1)]^{(\la)}$ contains both $(\0,h_0)$ and $w$. Thanks to \eqref{eq:apertquasiparab} and \eqref{eq:bounddist}, which are valid for $t \geq - h_0$,
we have $e_\la(v_1,\0)\le  c\sqrt{t}/R_\la$ and
$e_{\la}(v,v_1)\le c \sqrt{h_1-h}/ R_{\la} \le c \sqrt{t-h}/ R_{\la}$. Consequently, there exists a constant $c>0$ such that
\begin{equation}
  \label{eq:diamunion}
R_\la e_\la(v,\0)\le c(\sqrt{2t}+\sqrt{t-h}).
\end{equation}
There is a constant $c>0$ such that if $h \in [-ct^2, \infty)$, then $|v|=R_\la e_\la(v,\0)\le t$.

Consequently, when $\P^{(\la)}\cap {\mathcal U}\cap (\R^{d-1}\times (-\infty,-ct^2))=\emptyset$, then the localization radius of $\xi^{(\la)}$ is less than $t$. This means that
$$
 P[R^{\xi^{(\la)}}[(\0,h_0)] >  t] \le P[H((\0,h_0),\P^{(\la)})\ge t]+P[\P^{(\la)}\cap {\mathcal U}\cap (\R^{d-1}\times (-\infty,-ct^2))\ne\emptyset].
$$
Given \eqref{eq:diamunion}, we may use the same method as in \eqref{covsphere} to obtain
\begin{align*}
& d\P^{(\la)}(\P^{(\la)}\cap {\mathcal U}\cap (\R^{d-1}\times (-\infty,-c t^2))) \\& \le
c\int_{ct^2}^{\infty}e^{-ch}(\sqrt{t}+\sqrt{t+h})^{(d-1)}dh \le c e^{-t^2/c}. \ \ \  \qed
\end{align*}

\vskip.5cm

\noindent{\bf 4.2. Moment bounds for $\txi, \la \in [\la_0, \infty]$.}
 For
a random variable $X$ and $p \in (0, \infty)$, we let $||X||_p:= (\E
|X|^p)^{1/p}$. 

\begin{lemm}\label{L1}
For all $p \in [1, \infty)$ and $\xi\in \Xi$,
there is a constant $c>0$
such that for all $(v,h)  \in W_\la$, $\la \in [\la_0, \infty]$, we have
\begin{equation}\label{LIMITBD2}
      \E [| \xi^{(\la)}( (v,h),\P^{(\la)})|^p] \leq c  |h|^{c} \exp( - \frac{e^{h \vee 0}}{c} ).
\end{equation}
\end{lemm}

\noindent{\em Proof.}
We first prove \eqref{LIMITBD2} for a $k$-face functional $\xi^{(\la)}:= \xi_k^{(\la)}$, $k\in \{0,1,...,d-1\}.$
 We start by showing for all $\la \in [\la_0, \infty]$ and $h \in \R$
\begin{equation}\label{LIMITBD1-new}
   \sup_{v \in \R^{d-1}} \E [|\xi^{(\la)}((v,h),\P^{(\la)} )|^p] \leq c |h|^{c}.
\end{equation}
Since $\xi(x,\P_{\la})\eqd\xi(y,\P_{\la})$ whenever $|x|=|y|$, it follows that for all $(v,h)\in W_{\la}$, $$\xi^{(\lambda)}((\0,h),\P^{(\la)})\eqd \xi^{(\lambda)}((v,h),\P^{(\la)}).$$ Consequently, without loss of generality we may put  $(v,h)$ to be $(\0,h_0)$.

Let
$N^{(\la)}:=N^{(\la)}((\0,h_0)):= {\rm{card}}\{ {\rm{Ext}}^{(\la)}(\P^{(\la)} ) \cap  C(\0, R) \}$  with $R:=R^{ \xi^{(\la)} }[(\0,h_0)]$ the radius
of spatial localization for $\xi^{(\la)}$ at $(\0,h_0)$.  Clearly
$$
\xi^{(\la)}((\0,h_0),\P^{(\la)})\leq \frac{1}{k+1}\binom{N^{(\la)}((\0,h_0))}{k}.
$$
To show \eqref{LIMITBD1-new}, given $p \in [1, \infty)$, it suffices to show there is a constant $c:=c(p,k,d)$ such that for all $\la \in [\la_0, \infty]$
\begin{equation}\label{momentsN}
\E N^{(\la)}((\0,h_0))^{pk} \leq c  |h_0|^{c}.
\end{equation}

By \eqref{3d}, for all $r\in [0,\pi R_{\la}]$ and $\ell\in (-\infty,R_\la^2]$ we have
$$d\P^{(\la)}(C(\0,r)\cap (-\infty,\ell))\le cr^{d-1}(-\ell\vee 1)^{c}e^{\ell}.$$
Consequently, with $H:=H((\0,h_0),\P^{(\la)})$ as in Lemma \ref{BLem1} and ${\rm{Po}}(\alpha)$ denoting a Poisson random variable
with mean $\alpha$, we have for $\la \in [\la_0, \infty]$
\begin{align*}
& \ \ \ \  \E N^{(\la)}((\0,h_0))^{pk} \\
& \leq \E[ {\rm{card}}( \P^{(\la)} \cap [C(\0,R)\cap (-\infty,H)])^{pk}] \\
& = \sum_{i= 0}^{\infty}\sum_{j=h_0}^{\infty} \E[ {\rm{Po}}(d\P^{(\la)}(C(\0,R)\cap (-\infty,H)))^{pk} {\bf{1}}(i
\leq R < i + 1, j \leq H < j + 1) ] \\
& \leq \sum_{i= 0}^{\infty}\sum_{j=h_0}^{\infty}\E[ {\rm{Po}}(c(i + 1)^{d-1}(-(j+1)\vee 1)^ce^{(j+1)})^{pk} {\bf{1}}(R \geq  i, H \geq j) ].
\end{align*}

We shall repeatedly use the moment bounds for Poisson random variables, namely
$\E[ {\rm{Po}}(\alpha)^r] \leq c(r) \alpha^r, r \in [1, \infty).$ Using H\"older's inequality, we get
\begin{equation*}
\E N^{(\la)}((\0,h_0))^{pk} \le c\sum_{i= 0}^{\infty}\sum_{j=h_0}^{\infty}(i + 1)^{d-1}(-(j+1)\vee 1)^{{cpk}}e^{(j+1)pk}P[R\ge i]^{1/3}P[H\ge j]^{1/3}.
\end{equation*}
Splitting the sum on the $i$ indices into $i \in [0, |h_0|]$
and $i \in (|h_0|, \infty)$
yields with the help of Lemmas \ref{lem4.0} and \ref{BLem1}(a)
\begin{align*}
\E N^{(\la)}((\0,h_0))^{pk} & \leq  c |h_0|^c\sum_{j=0}^{\infty}e^{(j+1)pk}e^{-e^{j}/c}
+c\sum_{i = |h_0|}^{\infty} \sum_{j = 0}^{\infty}i^ce^{-i^2/c}e^{(j+1)pk}e^{-e^{j}/c}\\
&\le c|h_0|^c.
\end{align*}
This yields the required bound \eqref{momentsN}.

To deduce \eqref{LIMITBD2}, we argue as follows.
First consider the case $h_0 \in [0, \infty)$.
  By the Cauchy-Schwarz inequality and \eqref{LIMITBD1-new}
  \begin{align*}
  & \ \ \ \E  [|\xi^{(\la)}((\0,h_0),\P^{(\la)})|^p] \\
  & \leq
  (\E|\xi^{(\la)}((\0,h_0),\P^{(\la)})|^{2p})^{1/2} P[
  |\xi^{(\la)}((\0,h_0),\P^{(\la)})| > 0]^{1/2} \\
  & \leq (c(2p,k,d))^{1/2} |h_0|^{c_1(p, k, d)} P[
  |\xi^{(\la)}((\0,h_0),\P^{(\la)})| \neq 0]^{1/2}.
  \end{align*}
  The event $\{|\xi^{(\la)}((\0,h_0),v)| \neq
  0\}$ is a subset of the event that $(\0,h_0)$ is extreme in $\P^{(\la)}$ and
  we may now apply \eqref{fbd-1} for $t = h_0$, which is possible since we have
  assumed $h_0$ is positive.  This gives \eqref{LIMITBD2} for $h_0 \in [0, \infty)$.
  When $h_0 \in (-\infty, 0)$ we  bound $P[
  |\xi^{(\la)}((\0,h_0),\P^{(\la)})| > 0]^{1/2}$ by $c \exp( - e^0/c)$, $c$ large,
  which shows \eqref{LIMITBD2} for $h_0 \in (-\infty,0)$.
This concludes the proof of \eqref{LIMITBD2} when $\xi$ is a $k$-face functional.


We now prove \eqref{LIMITBD2} for $\xi_V$. For all $L \in (0, \infty)$ and $\la \in [\la_0, \infty)$, we put
$D^{(\la)}(L):=  ||\partial
\Phi^{(\la)}(\P^{(\la)}) \cap C(\0,L)||_{\infty}.$
Put $R:=R^{\xi_V^{(\la)}}[(\0,h_0)]$. 
The identity \eqref{Leb} shows that  $|\xi^{(\la)}_V((\0,h_0),\P^{(\la)})|$
is  bounded by the product of $c(1 + D^{(\la)}(R)/R_\la^2)^{d-1}$ and the Lebesgue measure of $B(\0,
R) \times [-D^{(\la)}(R), D^{(\la)}(R)]$.  We have
$$
\E | \xi^{(\la)}_V((\0,h_0),\P^{(\la)}) |^p \leq c \E(R^{d-1} D^{(\la)}(R)^{d})^p \leq c
|| R^{p(d-1)} ||_2 ||D^{(\la)}(R)^{pd}||_2,$$ by the Cauchy-Schwarz inequality.
By the tail behavior for $R$ we have  $\E R^r = r \int_0^{\infty} P[R > t] t^{r-1} dt \leq c(r)
|h_0|^r$ for all $r \in [1, \infty)$.  Also, for all $r \in [1, \infty)$ we have
$$
\E (D^{(\la)}(R))^r = \sum_{i = 0}^{\infty} \E (D^{(\la)}(R))^r {\bf{1}}(i \leq R < i
+ 1) \leq \sum_{i = 0}^{\infty} ||D^{(\la)}(i + 1)^r||_2 P[R \geq i]^{1/2}.
$$
By Lemma \ref{BLem1} we have $||D^{(\la)}(i + 1)^r||_2 \leq c(r) (i +
1)^{2(d-1)r}$, $\la \in [\la_0, \infty).$  We also have that $P[R \geq i], \ i \geq |h_0|,$ decays exponentially fast, showing that
$\E (D^{(\la)}(R))^r \leq c(r) |h_0|^{2(d-1)r}.$
It follows that
$$
\E | \xi^{(\la)}_V((\0,h_0),\P^{(\la)}) |^p \leq || R^{p(d-1)} ||_2 ||D^{(\la)}(R)^{pd}||_2
\leq c(p,d)|h_0|^{2pd(d-1)}  |h_0|^{(d-1)/2},
$$
which gives
\begin{equation}\label{LIMITBD3-new}
    \E |\xi^{(\la)}_V((\0,h_0),\P^{(\la)} )|^p \leq c |h_0|^{c}.
\end{equation}
The bound \eqref{LIMITBD2} for $\xi_V^{(\la)}$ follows from
\eqref{LIMITBD3-new} in the same way that \eqref{LIMITBD1-new}
implies \eqref{LIMITBD2} for $\xi_k^{(\la)}$.
 \qed

\vskip.5cm

\vskip.5cm

\noindent{\bf 4.3. Scaling limits.} The next two lemmas justify the assertion that functionals in $\Xi^{(\infty)}$ are indeed scaling limits of their counterparts in $\Xi^{(\la)}$.

\begin{lemm}  \label{L2a} For all  $h_0 \in \R$, $r \in (0, \infty)$, and $\xi\in\Xi$
we have  \be \label{rconv} \liml \E
\xi_{[r]}^{(\la)}((\0, h_0),\P^{(\la)}) =
 \E\xi^{(\infty)}_{[r]}((\0, h_0),\P).\ee
\end{lemm}

\noindent{\em Proof.}
Put $w_0:= (\0, h_0)$ and put $S(r,l):= B_{d-1}(\0,r) \times [-l, l]$, with $l$ a fixed deterministic height.
By Lemma \ref{BLem1} and the  Cauchy-Schwarz inequality, it is enough to show $$\liml \E
\xi_{[r]}^{(\la)}(w_0,\P^{(\la)}
 \cap S(r,l)) = \E\xi^{(\infty)}_{[r]}(w_0,\P \cap S(r,l)).$$
It is understood that the left-hand side is determined by the
geometry of the quasi-paraboloids $\{[ \Pi^{\uparrow}(w)]^{(\la)}
\}, w \in \P^{(\la)} \cap S(r,l),$ and similarly for the right-hand
side. Equip the collection $\X(r,l)$ of locally finite point sets in
$S(r,l)$ with the discrete topology. Thus if $\X_i, i \geq 1,$
is a sequence in $\X(r,l)$ and if
\be \label{seq} \lim_{i \to \infty}\X_i = \X, \ \ {\rm{then}} \ \ \X_i = \X \ {\rm{for}} \ \ i \geq i_0.
\ee
Recall that
$[\Pi^{\downarrow}(w)]^{(\infty)}$ coincides with
$\Pi^{\downarrow}(w).$  For all $\la \in [\la_0, \infty], w_1 \in W_\la,$ and $\X \in \X(r,l)$ we define $g_{k, \la}: W_\la \times
\X(r,l) \to \R$ by
taking $g_{k, \la}(w_1, \X)$ to be the product of $(k +
1)^{-1}$ and the number of quasi parabolic $k$-dimensional faces of
$\bigcup_{w \in \X} [\Pi^{\downarrow}(w)]^{(\la)}$ which contain
$w_1$, if $w_1$ is a vertex in  $\X$, otherwise $g_{k, \la}(w_1, \X)= 0$. Thus $g_{k, \la}(w_1, \X):= \xi_{[r]}^{(\la)}(w_1, \X \cap
S(r,l))$.

 Let $\X$ be in regular position, that is to say the intersection of $k$ quasi-paraboloids contains at most $(d-k+1)$ points
 of $\X$ for all $1\le k\le d$. Thus $\P$ is in regular position with probability one. To apply
 the continuous mapping theorem (Theorem 5.5 in \cite{Bi}),
  by \eqref{seq}, it is enough to show that $g_{k, \la}(w_0,\X)$ coincides with $g_{k, \infty}(w_0,\X)$ for $\lambda$ large enough. Let $\varepsilon>0$ be the minimal distance between any down paraboloid containing $d$ points of $\X$ and the rest of the point set.
 Perturbations of the paraboloids within an $\varepsilon$ parallel set do not change the number of  $k$-dimensional faces. In particular, for $\lambda$ large enough, the set $\partial\left(\cup_{w\in \X} [\Pi^{\downarrow}(w)]^{(\lambda)}\right)$ is included in that parallel set so that the number of $k$-dimensional faces does not change.
Thus $g_{k, \la}(w_0,\X)$ coincides with $g_{k, \infty}(w_0,\X)$ for large $\lambda$.


Since $\P^{(\la)} \tod \P$, we may apply the continuous mapping
theorem to get $$\xi_{[r]}^{(\la)}(w_0,\P^{(\la)}) \tod
 \xi^{(\infty)}_{[r]}(w_0,\P)$$ as $\la \to \infty$.
The convergence in distribution extends to convergence of expectations by the
uniform integrability of  $\xi_{[r]}^{(\la)}$, which follows from moment bounds for
$\xi_{[r]}^{(\la)}(w_0,\P^{(\la)})$ analogous to those for $\xi_{k}^{(\la)}(w_0,\P^{(\la)})$
as given in Lemma \ref{L1}.  This proves \eqref{rconv} when $\xi$ is a generic $k$-face functional.

Next we show for $\xi:= \xi_V$, $r \in (0, \infty)$ that
$$\liml \E \xi_{[r]}^{(\la)}(w_0,\P^{(\la)}
 \cap S(r,l)) = \E[\xi^{(\infty)}_{[r]}(w_0,\P \cap S(r,l)).$$
 This will yield \eqref{rconv}.
Recall that $\Vol_d^{(\la)}$ is the image of $R_\la \Vol_d$ under $T^{(\la)}$. Recall from \eqref{re-vol} the definition of
$\rm{Cyl}^{(\lambda)}(w)$.
For $\la \in [\la_0, \infty]$,  we define this time $\tilde{g}_{k, \la}: \R^{d-1} \times \R  \times
S(r,H) \mapsto \R$ by
\begin{align}   \label{yy}
&\hspace*{-.3cm}\tilde{g}_{k,\la} (w,\X) = \xi_{[r]}^{(\la)}(w,\X
 \cap S(r,l))\nonumber\\
&\hspace*{-.3cm} =  {\rm{Vol}}_d^{(\la)}
(\{(v,h) \in S(r,l): 0 \leq h \le \partial \Phi^{(\la)}(\X)(v), v\in {\rm{Cyl}}^{(\lambda)}(w),
\Phi^{(\la)}({\X})(v) \geq 0\})\nonumber\\
&\hspace*{-.3cm}
- {\rm{Vol}}_d^{(\la)}
(\{(v,h) \in S(r,l): \Phi^{(\la)}(\X)(v) \leq h \le 0, v\in {\rm{Cyl}}^{(\lambda)}(w),
\Phi^{(\la)}(\X)(v) < 0\}).
\end{align}
Recalling \eqref{seq}, it is enough to show for a fixed point set $\X$ in regular position that
$$
\liml | \tilde{g}_{k, \la}(w,\X) - \tilde{g}_{k, \infty}(w,\X) | = 0.
$$
We show that the first term in \eqref{yy} comprising $\tilde{g}_{k, \la}(w,\X)$ converges to the
first term comprising $\tilde{g}_{k, \infty}(w,\X)$. In other words,
setting for all $\la \in [\la_0, \infty)$
$$
F^{(\la)}(\X) := \{(v,h) \in S(r,l): 0 \leq h \le \partial \Phi^{(\la)}(\X)(v), v\in {\rm{Cyl}}^{(\lambda)}(w),
\Phi^{(\la)}(\X)(v) \geq 0\}$$
and writing $F(\X)$ for $F^{(\infty)}(\X)$, we show
$$
\liml | {\rm{Vol}}_d^{(\la)}
(F^{(\la)}(\X)) - {\rm{Vol}}_d
(F(\X)) | = 0.$$
The proof that the second  term comprising $\tilde{g}_{k, \la}(w,\X)$ converges to the
second term comprising $\tilde{g}_{k, \infty}(w,\X)$ is identical.
We have
\be \label{T1}
| {\rm{Vol}}_d^{(\la)}
(F^{(\la)}(\X)) - {\rm{Vol}}_d
(F(\X)) | \leq | {\rm{Vol}}_d^{(\la)}
(F^{(\la)}(\X)) - {\rm{Vol}}_d^{(\la)}
(F(\X)) | \ee
$$ + \ | {\rm{Vol}}_d^{(\la)}
(F(\X)) - {\rm{Vol}}_d
(F(\X)) |.
$$
Since $\partial \Phi^{(\la)}(\X)$ converges uniformly to
$\partial \Phi(\X)$ on compacts (recall Lemma \ref{paralem}; see also the proof of Proposition \ref{BLem2} below) and since $d^H({\rm{Cyl}}^{(\lambda)}(w), {\rm{Cyl}}(w))$ decreases to zero as $\la \to \infty$ (indeed $\partial({\rm{Cyl}}^{(\la)}(w)) \to \partial {\rm{Cyl}}(w)$ uniformly), we get
for $\la \in [\la_0, \infty)$ that $F^{(\la)}(\X) \Delta F(\X)$ is a subset of a set
$A(\X) \subset \R^d$ of arbitrarily small volume.  So
$| {\rm{Vol}}_d^{(\la)}
(F^{(\la)}(\X) ) - {\rm{Vol}}_d^{(\la)}
(F(\X)) | \leq {\rm{Vol}}_d^{(\la)} (A(\X)).$
By Lemma \ref{weakcon}, we have  $\rm{Vol}_d^{(\la)} \tod \rm{Vol}_d$ and thus the first term in
\eqref{T1}  goes to zero as $\la \to \infty$.  Appealing again to $\rm{Vol}_d^{(\la)} \tod \rm{Vol}_d$,
 the second term in \eqref{T1} likewise tends to zero, showing \eqref{rconv} as desired.
\qed

\vskip.5cm

\begin{lemm}  \label{L2} For all $h_0 \in \R$ and $\xi\in \Xi$
we have
$$\liml \E\xi^{(\la)}((\0,h_0),\P^{(\la)})=
 \E\xi^{(\infty)}((\0,h_0),\P).$$
\end{lemm}

 \noindent{\em Proof.}    Let $w_0:=(\0,h_0)$.   By Lemma \ref{L2a}, given $\epsilon > 0$, we have
 for all $\la  \in [\la_0(\epsilon), \infty)$
 \be \label{tri1}
 |\E \xi_{[r]}^{(\la)}(w_0, \P^{(\la)}) - \E \xi_{[r]}^{(\infty)}(w_0, \P))| < \epsilon.
 \ee
We now show that replacing $\xi_{[r]}^{(\la)}$ and $\xi_{[r]}^{(\infty)}$ by
$\xi^{(\la)}$ and $\xi^{(\infty)}$, respectively, introduces negligible error
in \eqref{tri1}.
 Write
\begin{align*}
& \ \ \ \ |\E \xi_{[r]}^{(\la)}(w_0, \P^{(\la)}) - \E \xi^{(\la)}(w_0, \P^{(\la)}))| \\
& = |\E (\xi_{[r]}^{(\la)}(w_0, \P^{(\la)}) - \xi^{(\la)}(w_0, \P^{(\la)})){\bf 1}(R^{\xi^{(\la)}}[w_0] < r)| \\
&  \ \ \ + \ |\E (\xi_{[r]}^{(\la)}(w_0, \P^{(\la)}) - \xi^{(\la)}(w_0, \P^{(\la)})){\bf 1}(R^{\xi^{(\la)}}[w_0] > r)|.
 \end{align*}
 The first term vanishes by definition of $R^{\xi^{(\la)}}[w_0]$. By the Cauchy-Schwarz inequality and  Lemma \ref{lem4.01}, the second term is bounded by
 \be \label{CS1}
 || \xi_{[r]}^{(\la)}(w_0, \P^{(\la)}) - \xi^{(\la)}(w_0, \P^{(\la)}) ||_2 P[R^{\xi^{(\la)}}[w_0] > r]^{1/2}
 \leq c P[R^{\xi^{(\la)}}[w_0] > r]^{1/2} \leq \epsilon
 \ee
 if $r \in [|h_0|, \infty)$ is large enough.
 For $r \in [r_0(\epsilon, h_0), \infty)$ and $\la \in [\la_0(\epsilon),\infty)$ it follows that
 \be \label{tri2}
 |\E \xi^{(\la)}(w_0, \P^{(\la)}) - \E \xi^{(\la)}_{[r]}(w_0, \P^{(\la)})| < \epsilon.
 \ee
 Similarly for  $r \in [r_1(\epsilon, h_0), \infty)$ we have
 \be \label{tri3}
 |\E \xi^{(\infty)}(w_0, \P) - \E \xi^{(\infty)}_{[r]}(w_0, \P)| <  \epsilon.
 \ee
 Combining \eqref{tri1}-\eqref{tri3} and using the triangle inequality we get
 for $r \geq (r_0(\epsilon) \vee r_1(\epsilon))$ and $\la \in [\la_0(\epsilon), \infty)$
 $$
 |\E \xi^{(\la)}(w_0, \P^{(\la)}) - \E \xi^{(\infty)}(w_0, \P)| < 3 \epsilon.
 $$
Lemma \ref{L2} follows since $\epsilon$ is arbitrary.    \qed

\vskip.5cm

\noindent{\bf 4.4. Two point correlation function for $\txi$.} For
all $h \in \R$, $(v_1,h_1) \in W_\la,$ and $\xi\in\Xi$
we extend definition \eqref{SO2} by putting for
all $\la \in [\la_0, \infty]$
$$
c^{\xi^{(\la)}}((\0,h_0), (v_1,h_1)) :=
$$
$$
\E[ \xi^{(\la)}((\0,h_0),\P^{(\la)} \cup \{(v_1,h_1)\}) \
\xi^{(\la)}((v_1,h_1),\P^{(\la)} \cup \{(\0, h_0)\} )] - $$
$$\E \xi^{(\la)}((\0,h_0),\P^{(\la)} ) \E
\xi^{(\la)}((v_1,h_1),\P^{(\la)}).
$$

The next lemma shows convergence of the re-scaled  two-point correlation functions on
re-scaled input $\P^{(\la)}$ to their counterpart correlation functions on the limit input $\P$.

\begin{lemm}  \label{L2-new} For all $h_0 \in \R$, $(v_1,h_1) \in \R^{d-1} \times \R,$ and $\xi \in \Xi$  we have
$$\liml c^{\xi^{(\la)}}((\0,h_0), (v_1,h_1) ) =
c^{\xi^{(\infty)}}((\0,h_0), (v_1,h_1) ).$$
\end{lemm}

\noindent{\em Proof.} 
We deduce from Lemma \ref{L2} that
$$\liml \E \xi^{(\la)}((\0,h_0),\P^{(\la)} ) \E
\xi^{(\la)}((v_1,h_1) ,\P^{(\la)})=\E \xi^{(\infty)}((\0,h_0),\P) \E
\xi^{(\infty)}((v_1,h_1) ,\P).$$
By the Cauchy-Schwarz inequality, we get
\begin{align*}
&| \E[ \xi^{(\la)}((\0,h_0),\P^{(\la)} \cup \{(v_1,h_1)\} ) \
\xi^{(\la)}((v_1,h_1) ,\P^{(\la)} \cup \{(\0, h_0)\})]\\&\hspace*{1cm}-\E[ \xi^{(\infty)}((\0,h_0),\P \cup \{(v_1,h_1)\} ) \ \
\xi^{(\infty)}((v_1,h_1) ,\P \cup \{(\0, h)\})]|\le T_1(\la) +T_2(\la),
\end{align*}
where
\begin{align*}
T_1(\la)&:=\E[|\xi^{(\la)}((\0,h_0),\P^{(\la)} \cup \{(v_1,h_1)\} )-\xi^{(\infty)}((\0,h_0),\P \cup \{(v_1,h_1)\} )|^2]^{1/2} \\ & \ \ \times \E[|\xi^{(\la)}((v_1,h_1) ,\P^{(\la)} \cup \{(\0, h_0)\})|^2]^{1/2}
\end{align*}
and
\begin{align*}
T_2(\la) &:=\E[|\xi^{(\la)}((v_1,h_1) ,\P^{(\la)} \cup \{(\0,h_0)\})-\xi^{(\infty)}((v_1,h_1) ,\P \cup \{(\0,h_0)\})|^2]^{1/2}\\ & \ \ \times \E[|\xi^{(\infty)}((\0,h_0),\P^{(\la)} \cup \{(v_1, h_1)\})|^2]^{1/2}.
\end{align*}
Throughout we let $\P^{(\la)}, \la \geq 1,$ and $\P$
be defined on the same probability space, with $\P^{(\la)}$ independent of $\P$
for all $\la \geq 1$.  We couple $\P^{(\la)}$ and $\P$ so that
\begin{align}
  \label{eq:jointdist}
&\hspace*{-1cm}(\xi^{(\la)}((v_1,h_1) ,\P^{(\la)} \cup \{(\0,h_0)\}), \xi^{(\infty)}((v_1,h_1) ,\P \cup \{(\0, h_0)\})
\nonumber\\&\overset{D}=(\xi^{(\la)}((\0,h_0),\P^{(\la)} \cup \{(-v_1,h_0)\}), \xi^{(\infty)}((\0,h_0),\P \cup \{(-v_1, h_0)\}).
\end{align}
We  show first  $\liml T_1(\la)= 0$.
We have seen in the proof of Lemma \ref{L2a} that $\xi_{[r]}^{(\la)}((\0,h_0),\P^{(\la)} \cup \{(v_1,h_1)\} )$ converges in distribution to $\xi_{[r]}^{(\infty)}((\0,h_0),\P \cup \{(v_1,h_1)\} )$ for every $r>0$. Lemma \ref{L1} implies that this family is uniformly integrable so the convergence occurs in $L^2$, that is to say
$$
\liml (\E | \xi_{[r]}^{(\la)}((\0,h_0),\P^{(\la)} \cup \{(v_1,h_1)\} ) - \xi_{[r]}^{(\infty)}((\0,h_0),\P \cup \{(v_1,h_1)\} )|^2)^{1/2} = 0.
$$
 Using the same method as in the proof of Lemma \ref{L2}, we obtain
\begin{equation}
  \label{eq:lemma4.6}
\liml \E[|\xi^{(\la)}((\0,h_0),\P^{(\la)} \cup \{(v_1,h_1)\} )-\xi^{(\infty)}((\0,h_0),\P \cup \{(v_1,h_1)\} )|^2]=0.
\end{equation}
By Lemma \ref{L1}, the variables $\xi^{(\la)}((v_1,h_1) ,\P^{(\la)} \cup \{(\0, h_0)\})$ are uniformly bounded in $L^2$ so we deduce from \eqref{eq:lemma4.6} that  $\liml T_1(\la)= 0$.
To see that $\liml T_2(\la)= 0$, we use \eqref{eq:jointdist} and we follow the proof that $\liml T_1(\la)= 0$.
\qed

\vskip.5cm

The next lemma shows that the re-scaled and limit two point correlation function decays exponentially fast with the distance between spatial coordinates of the input and
super-exponentially fast with respect to positive height coordinates.

\begin{lemm} \label{L3} For all $\xi \in \Xi$
there is a constant $c:=c(\xi,d) \in (0, \infty)$ such that
for all  $(v_1,h_1) \in W_{\la}$ satisfying  $|v_1| \geq 2 \max( |h_0|, |h_1|)$ and $\la \in [\la_0, \infty]$ we have
\be \label{4.7a} |c^{\xi^{(\la)}}((\0,h_0), (v_1,h_1))| \leq c |h_0|^{c}
|h_1|^{c_3}
\exp\left( {-1 \over c} (|v_1|^2+e^{h_0 \vee 0}+e^{h_1 \vee 0}) \right).
\ee
\end{lemm}

\noindent{\em Proof.}
Let $x_\la:= [T^{(\la)}]^{-1}((\0,h_0))$ and $y_\la:= [T^{(\la)}]^{-1}((v_1,h_1))$.
Put
$$X_\la:= \xi^{(\la)}((\0,h_0),\P^{(\la)} \cup (v_1,h_1))= \xi(x_\la, \P_\la \cup y_\la),$$
$$Y_\la:= \xi^{(\la)}((v_1,h_1),\P^{(\la)} \cup (\0,h))= \xi(y_\la, \P_\la \cup x_\la),$$
$$\tX_\la:= \xi^{(\la)}((\0,h_0),\P^{(\la)})= \xi(x_\la, \P_\la),$$
$${\rm{and}} \  \tY_\la:= \xi^{(\la)}((v_1,h_1),\P^{(\la)})= \xi(y_\la, \P_\la).$$
We have
$$
|c^{\xi^{(\la)}}((\0,h_0), (v_1,h_1))| = |\E X_\la Y_\la - \E {\tX}_\la \E {\tY}_\la|
$$
which gives for all $r \in (0, \infty)$
\begin{align*} |c^{\xi^{(\la)}}((\0,h_0), (v_1,h_1))|
& \leq |\E X_\la Y_\la {\bf{1}} ( R^{\xi}(x_\la, \P_\la) \leq r,  R^{\xi}(y_\la, \P_\la) \leq r) \\
& - \E \tX_\la {\bf{1}} (  R^{\xi}(x_\la, \P_\la) \leq r) \E \tY_\la {\bf{1}}( R^{\xi}(y_\la, \P_\la) \leq r)| \\
& + |\E X_\la Y_\la [{\bf{1}} [ ( R^{\xi}(x_\la, \P_\la) \geq r) + {\bf{1}}  (R^{\xi}(y_\la, \P_\la) \geq r)]| \\
& + |\E \tX_\la  \E \tY_\la {\bf{1}}( R^{\xi}(x_\la, \P_\la) \geq r)|  +
|\E \tX_\la  \E \tY_\la {\bf{1}} ( R^{\xi}(y_\la, \P_\la) \geq r)|.
\end{align*}
 Put  $r := |v_1|/2R_\la$.
This choice of $r$ ensures that the difference of the first two terms is zero
by independence of $X_\la {\bf{1}} ( R^{\xi}(x_\la, \P_\la) \leq r)$
and $ Y_\la {\bf{1}} (  R^{\xi}(y_\la, \P_\la) \leq r)$.
Recall that  $R^{\xi}(x_\la, \P_\la)$ and  $R^{\xi}(y_\la, \P_\la)$ have the same distribution.  When $|v_1| \geq 2 \max( |h_0|, |h_1|)$, H\"older's inequality implies that the third term is bounded by
\begin{align}
&||X_\la||_3 ||Y_\la||_3 [ P [R^{\xi}(x_\la, \P_\la) \geq r]^{1/3} +
P [R^{\xi}(y_\la, \P_\la) \geq r]^{1/3} ]\nonumber\\
& \le c |h_0|^{c}
|h_1|^c
\exp\left( {-1 \over c} (e^{h_0 \vee 0}+e^{h_1 \vee 0}) \right)P [R^{\xi}(x_\la, \P_\la) \geq r]^{1/3} \label{eq:domin}\\
&\leq c |h_0|^{c}
|h_1|^c
\exp\left( {-1 \over c} (|v_1|^2+e^{h_0 \vee 0}+e^{h_1 \vee 0})\right).\nonumber
\end{align}
The fourth and fifth terms are bounded similarly, giving \eqref{4.7a}.  \qed

\section{Proofs of main results}\label{sectionproofs}

\allco

\noindent{\bf 5.1. Proof of Theorems  \ref{Th0} and \ref{Th4}.}
The next result contains  Theorem \ref{Th4} and it yields Theorem \ref{Th0},
since it implies that the set ${\rm{Ext}}^{(\la)}(\P^{(\la)})$ of extreme points of $\P^{(\la)}$ converges in law to  ${\rm{Ext}}(\P)$ as $\la \to \infty$ (indeed, the set ${\rm{Ext}}^{(\la)}(\P^{(\la)})$ is also the set of local minima of the function $\partial \Phi^{(\la)}(\P^{(\la)})$).

\begin{prop}\label{BLem2} Fix $L \in (0, \infty).$  The boundary  of $\Psi^{(\la)}(\P^{(\la)})$
converges in probability as $\la \to \infty$ to the boundary  of $\Psi(\P)$
in the space $\C(B_{d-1}(\0,L)))$ equipped with the supremum norm.
Similarly, the boundary of $\Phi^{(\la)}(\P^{(\la)})$ converges in probability  as
$\la \to \infty$ to the Burgers' festoon $\partial(\Phi(\P))$.
\end{prop}

\noindent{\em Proof.} We only prove the first convergence statement as the second  is
handled similarly.  We show for fixed $L \in (0, \infty)$ that the
boundary of $\Psi^{(\la)}(\P^{(\la)})$ converges in law to $
\partial( \Psi(\P))$ in the space $\C(B_{d-1}(\0,L))$. With $L$
fixed, for all $l \in [0, \infty)$ and $\la \in [0, \infty)$, let
$E(L, l, \la)$ be the event that the heights of
$\partial(\Psi^{(\la)}(\P^{(\la)}))$ and $\partial( \Psi(\P))$
belong to $[-l, l]$ over the spatial region $B_{d-1}(\0,L)$. By
Lemma \ref{BLem1}, we have that $P[E(L, l, \la)^c]$ decays
exponentially fast in $l$, uniformly in $\la$, and so it is enough
to show, conditional on $E(L, l, \la)$, that
$\partial(\Psi^{(\la)}(\P^{(\la)}))$ is close to
$\partial(\Psi(\P))$ in the space $\C(B_{d-1}(\0,L))$, $\la$ large.

Recalling the definition of $\Psi^{(\la)}(\P^{(\la)})$ at
\eqref{Boolmod}, we need to show, conditional on $E(L, l, \la)$,
that the boundary of
$$
\bigcup_{w \in \P^{(\la)} \cap C(\0, L) } (
[\Pi^{\uparrow}(w)]^{(\la)} \cap C(\0, L) )
$$
is close to the boundary of \be \label{bdry}
 \bigcup_{ \v \in \P \cap C(\0, L) } ( \Pi^{\uparrow}(\v) \cap C(\0, L) ).
\ee

By Lemma \ref{paralem}, given $w_1:=(v_1,h_1) \in \P^{(\la)} \cap C(\0, L)$, it follows that on $E(L, l, \la)$ the boundary of
$[\Pi^{\uparrow}(w_1)]^{(\la)} \cap C(\0, L)$ 
is within $O(R_\la^{-1})$ of the boundary of
$\Pi^{\uparrow}(w_1) \cap C(\0, L)$. The boundary of $\Psi^{(\la)}(\P^{(\la)}) \cap
C(\0, L)$ is a.s. the finite union of graphs of the above form and
is thus a.s. within $O(R_\la^{-1})$ of the boundary of
$$
 \bigcup_{ \v \in \P^{(\la)} \cap C(\0, L) } ( \Pi^{\uparrow}(\v) \cap C(\0, L) ).
 $$
It therefore suffices to show that the boundary of $\bigcup_{ \v \in
\P^{(\la)} \cap C(\0, L) } ( \Pi^{\uparrow}(\v) \cap C(\0, L) )$ is
close to the boundary of the set given at \eqref{bdry}.  However, we may couple
$\P^{(\la)}$ and $\P$ on $B_{d-1}(\0, L) \times [-l, l]$ so that
they coincide except on a set with probability less than $\epsilon$,
showing the desired closeness with probability at least $1 -
\epsilon$.  \qed


\vskip.5cm

\noindent{\bf 5.2. Proof of expectation asymptotics \eqref{main1}.}
For $g \in {\cal C}(\S^{d-1})$ and $\xi \in \Xi$ we have
\be \label{dis1}
  \E  [\langle g_{R_\la}, \mu^{\xi}_{\la} \rangle] =  \int_{\R^d} g(  \frac{x} { R_\la } ) {\Bbb E}\left[\xi(
  x,\P_{\la}) \right]  \la \phi(x) dx.
 \ee

Since $\xi(x,\P_\la)\overset{D}{=}\xi(y,\P_\la)$ as soon as $|x|=|y|$, we have
$$\E\xi(x, \P_\la) = \E\xi(|x|u_0,\P_\la)=\E[\xi^{(\la)}\left(
   (\0,h_0),\P^{(\la)} \right)]$$
where $h_0$ is defined by $|x|=R_\la(1-h_0/R_\la^2)$. Writing $u=x/|x|$, we have by \eqref{3b}
$dx =
[ R_{\la}(1-h_0/R_{\la}^2)]^{d-1} R_\la^{-1} dh_0 d\sigma_{d-1}(u)$. Consequently,
we see that $R_\la^{-(d-1)}\E  [\langle g_{R_\la}, \mu^{\xi}_{\la} \rangle]$ transforms to
\begin{align*}
\int_{u \in  \S^{d-1} } \int_{h_0 \in (-\infty, R_\la^2]} g( u  (1 -
{h_0 \over R_\la^2})) {\Bbb E}\left[\xi^{(\la)}\left(
  (\0,h_0),\P^{(\la)} \right)\right] \tp_\la(u,h_0) (1 - {h_0 \over R_\la^2})^{d-1}  dh_0
d\sigma_{d-1}(u),
\end{align*}
where $\tp: \S^{d-1} \times \R \to \R$ is given by
\be \label{defG} \tp_\la(u,h_0) := \frac{\la} {R_\la}
\phi(u \cdot R_\la (1 - {h_0 \over R^2_\la} )).\ee

By \eqref{3c} we have for all $u \in \S^{d-1}$ that
$$
\tp_\la(u,h_0) = \frac{ \sqrt{2 \log \la} } {R_\la}  \exp\left(h_0 - \frac{h_0^2}{2R_{\la}^2}\right).
$$

Thus there is $c \in (0, \infty)$ such that for all $h_0 \in \R$ we have
\be \label{5.2}
  \sup_{u \in \S^{d-1} } \sup_{\la \geq 3}  \tp_\la(u,h_0) \leq ce^{h_0}
   \ee
   and for all $u \in \S^{d-1}$
   \be \label{5.3}
   \lim_{\la \to \infty} \tp_\la(u,h_0) = e^{h_0}.
   \ee

\vskip.3cm
By the continuity of $g$,  Lemma \ref{L2}  and the limit \eqref{5.3}, we have for  $h_0 \in (-\infty, R_\la^2)$ that the integrand
inside the double integral converges to $g(u) \E
\left[\xi^{(\infty)} \left(
  (\0,h_0),\P \right)\right] e^{h_0}$ as $\la \to \infty$.  Moreover,
by
  \eqref{5.2} and the moment bounds of Lemma \ref{L1}, the integrand is dominated by the product of a polynomial in $h_0$ and
an exponentially decaying function of $h_0$.
The dominated  convergence theorem gives the claimed result \eqref{main1}.  \qed

\vskip.5cm

\noindent{\bf 5.3. Proof of variance  asymptotics \eqref{main2}.}
For $g \in {\cal C}(\S^{d-1}),$ using the Mecke-Slivnyak formula (Corollary 3.2.3 in \cite{SW}), we have \be
 \Var  [\langle g_{R_\la}, \mu^{\xi}_{\la} \rangle] :=  I_1(\la) + I_2(\la), \label{dis3}
 \ee
 where
 $$
  I_1(\la) :=  \int_{\R^d} g(  \frac{x} { R_\la } )^2 {\Bbb E}\left[\xi(
  x,\P_{\la})^2 \right]  \la \phi(x) dx $$
  and
  $$
 I_2(\la):=
  \int_{\R^d} \int_{\R^d} g(  \frac{x} { R_\la } ) g(  \frac{y} { R_\la } )
[ \E \xi(x, \P_\la \cup {y} )\xi(y, \P_\la \cup {x} ) - \E \xi(x,
\P_\la )\E \xi(y, \P_\la ) ] \la^2 \phi(y)  \phi(x) dy dx.
$$

We examine $\lim_{\la \to \infty} R_\la^{-(d-1)} I_1(\la)$ and
$\lim_{\la \to \infty} R_\la^{-(d-1)} I_2(\la)$ separately.
As in the proof of expectation asymptotics \eqref{main1}, we have
\be \label{Var1} \lim_{\la \to \infty} R_\la^{-(d-1)} I_1(\la) =
\int_{-\infty}^\infty \E \xi^{(\infty)}_{k}((\0,h_0), \P)^2 e^{h_0} dh_0
\int_{\S^{d-1}} g(u)^2 du. \ee
Next consider $\lim_{\la \to \infty} R_\la^{-(d-1)} I_2(\la)$.  For $x \in \R^d$ we write
\be \label{param1}
x = u R_\la (1 - {h_0\over R_\la^2} ), \ (u,h_0) \in \S^{d-1} \times \R. 
\ee
We now re-scale the integrand in $I_2(\la)$ as follows.
Given $u:=u_x \in \S^{d-1}$ in the definition of $x$, define
$T^{(\la)}$  as in \eqref{scaltrans}, but with $u_0$ there replaced by $u$.
 Write
$T_u^{(\la)}$ to denote the dependency on $u$. Denoting by $(\0,h_0)$ and $(v_1,h_1)$ the images under $T_u^{(\la)}(x)$ of $x$ and $y$ respectively, we notice that $R_\la^{-(d-1)} I_2(\la)$ is transformed as follows.

\vskip.5cm

\noindent (i) The `covariance' term $[ \E \xi(x, \P_\la \cup {y} )\xi(y, \P_\la \cup {x} ) - \E \xi(x, \P_\la )\E \xi(y, \P_\la ) ]$ transforms to $ c^{\xi^{(\la)}}((\0,h_0),
(v_1,h_1))$.  By Lemma \ref{L2-new} we have uniformly in $v_1 \in
T_u^{(\la)}(\S^{d-1})$ and $h_0, h_1 \in \R$ that \be \label{lim3} \lim_{\la \to \infty}c^{\xi^{(\la)}}((\0,h_0),
(v_1,h_1)) =  c^{ \xi^{( \infty)} }((\0,h_0), (v_1,h_1)). \ee

\noindent (ii) The product  $g(  \frac{x} { R_\la } ) g(  \frac{y} { R_\la })$ becomes
$$f_{1, \la}(u, h_0, v_1,h_1):= g(u  (1 - {h_0 \over R_\la^2} )) g(R_\la^{-1}[T_u^{(\la)}]^{-1} ((v_1,h_1)) ).$$
Using \eqref{scaltrans} and \eqref{param1}, we notice that
 $[T_u^{(\la)}]^{-1} ((v_1,h_1))=R_\la\left(1-\frac{h_1}{R_{\la}^2}\right)\exp_{d-1}R_{\la}^{-1}v_1$ and consequently $$\lim_{\la\to\infty}R_\la^{-1}[T_u^{(\la)}]^{-1} ((v_1,h_1)) =u.$$
By continuity of $g$, we then have
uniformly in $v_1 \in T_u^{(\la)}(\S^{d-1})$ and $h_0, h_1 \in \R$ that \be \label{lim1} \lim_{\la \to
\infty}f_{1, \la}(u, h_0, v_1,h_1) = g(u)^2. \ee

\vskip.3cm

\noindent (iii) The double integral over $(x,y) \in \R^d \times \R^d$
transforms into a quadruple integral over $(u, h_0, v_1, h_1) \in
\S^{d-1} \times (-\infty, R_\la^2] \times T_u^{(\la)}(\S^{d-1}) \times
(-\infty, R_\la^2]$.

\vskip.3cm

\noindent (iv) By \eqref{3d}, the differential $\la \phi(y)dy$
transforms to
$$\frac{ \sin^{d-2}( R_\la^{-1} |v_1| )}{ |R_\la^{-1} v_1|^{d-2} }
{ \sqrt{2 \log \la} \over R_\la } (1 - {h_1 \over R_\la^2} )^{d-1}
\exp\left( h_1 - \frac{{h_1}^2}{2R_{\la}^2}\right)dv_1dh_1$$ whereas  $R_\la^{-(d-1)} \la \phi(x)dx$ transforms
to
$$
\tp_\la(u,h_0) (1 - {h_0 \over R_\la^2})^{d-1}  dh_0 d\sigma_{d-1}(u).
$$
Thus the product  $R_\la^{-(d-1)} \la^2 \phi(y)  \phi(x) dy dx$
transforms to
$$
f_{2, \la}(u,h_0, v_1,h_1)d\sigma_{d-1}(u) dh_0dv_1 dh_1
$$
where
$$
f_{2, \la}(u,h_0, v_1,h_1) := \frac{ \sin^{d-2}( R_\la^{-1} |v_1| )}{ |R_\la^{-1} v_1|^{d-2} }{ \sqrt{2 \log \la} \over R_\la } (1 - {h_1
\over R_\la^2} )^{d-1} \exp(h_1 - \frac{h_1^2}{2R_{\la}^2}) \tp_\la(u,h_0) (1 - {h_0 \over R_\la^2})^{d-1}.
$$
As in the proof of Lemma \ref{weakcon} and by \eqref{5.3} we have uniformly in $u \in
\S^{d-1}, v_1 \in T^{(\la)}(\S^{d-1})$ and $h_0, h_1 \in \R$ that \be \label{lim2} \lim_{\la \to
\infty}f_{2, \la}(u,h_0, v_1,h_1) = e^{h_0 + h_1}. \ee

\vskip.3cm

We re-write $R_\la^{-(d-1)} I_2(\la)$ as
 $$\displaystyle {= \int_{u \in  \S^{d-1}
} \int_{h_0 \in (-\infty, R_\la^2]} \int_{ T_u^{(\la)}(\S^{d-1}) }
\int_{h_1 \in (-\infty, R_\la^2]}F_\la(u,h_0, v_1,h_1)dh_1 dv_1 dh_0
d\sigma_{d-1}(u),}
$$
where
$$
F_\la(u, h_0, v_1,h_1) := f_{1, \la}(u,h_0, v_1,h_1) c^{\xi^{(\la)}}((\0,h_0), (v_1,h_1))
f_{2, \la}(u,h_0, v_1,h_1).
$$

Combining  the limits \eqref{lim3}-\eqref{lim2}, it follows for all
$u \in \S^{d-1}, h_0, h_1 \in \R$ and $v_1 \in T^{(\la)}(\S^{d-1})$,
that we have the pointwise convergence
$$
\lim_{\la \to \infty}F_\la(u,h_0, v_1,h_1) = g(u)^2 c^{ \xi^{( \infty)} }((\0,h_0), (v_1,h_1)) e^{h_0 + h_1}.
$$
By Lemma \ref{L3}, \eqref{eq:domin} and \eqref{5.2}, we get that
\begin{align*}
&|F_\la(u,h_0,v_1,h_1)|\le c |h_0|^{c}
|h_1|^{c}
\exp\left( {-1 \over c} (e^{h_0 \vee 0}+e^{h_1 \vee 0})+h_0 +h_1 \right)\\&\hspace*{4cm}\cdot\left(\exp({-1 \over c}|v_1|^2)+{\bf 1}(|v_1|\le 2\max(|h_0|,|h_1|))\right).
\end{align*}
Using that there exists $c'>0$ such that $\frac{1}{c}e^{h\vee 0}-h\ge \frac{1}{c'}e^{h\vee 0}$ holds for all $h \in \R$, we obtain that
  $F_\la(u, h_0, v_1,h_1)$ is dominated by an exponentially decaying
function of all arguments and is therefore integrable.
The dominated convergence theorem gives
\be \label{Var2}
 \lim_{\la \to \infty} R_\la^{-(d-1)} I_2(\la) = \ee
 $$
\int_{\S^{d-1} } \int_{\R^{d-1}}  \int_{h_0 \in (- \infty, \infty )}
\int_{h_1 \in (- \infty, \infty )} g(u)^2 c^{ \xi^{( \infty)} }((\0,h_0), (v_1,h_1)) e^{h_0 + h_1} dh_0 dh_1 dv_1  d\sigma_{d-1}(u). $$
Combining \eqref{Var1} and \eqref{Var2} gives the claimed variance
asymptotics \eqref{main2}. The positivity of $F_{k,d}$ is
established in \cite{BVu}, concluding the proof of \eqref{main2}.
  \qed

\vskip.5cm

\noindent{\bf 5.4. Proof of Corollary \ref{Th3}.}
Define
$\xi(x, \P)$ to be one if $x \in {\rm{Ext}}(P)$, otherwise
put $\xi(x, \P)= 0.$ Put
$$
\mu_\la:= \sum_{x \in {\rm{Ext}}(P) \cap Q_\la} \xi(x, \P \cap Q_\la) \delta_x.
$$
Note that $\E [ {\rm{card}} ({\rm{Ext}}( \P \cap  Q_\la))] = \E [\langle 1, \mu_\la \rangle].$ Put $\tQ_\la:= [ {-1 \over 2} \la^{1/(d-1)}, {1 \over 2} \la^{1/(d-1)}]^{d-1}$ so that $Q_\la = \tQ_\la \times \R$. Writing generic points in $\R^{d-1} \times \R$ as $(v,h)$ we have
$$
\la^{-1} \E [\langle 1, \mu_\la \rangle] = \la^{-1} \int_{\tQ_\la} \int_{- \infty}^{\infty}
 \E \xi( (v,h), \P \cap Q_\la) e^h dh dv.
 $$

  Put $\tP_\la$ to be a Poisson point process on $\R^{d-1} \times \R$ with intensity $\la e^h dh dv.$
  For any $\X \subset \R^{d-1} \times \R$, let $\alpha \X:= \{ (\alpha v, h): (v,h) \in \X \}.$
  Thus $\P \cap Q_\la \eqd \la^{1/(d-1)}[ \tP_\la \cap Q_1 ]$ and
  $$
  \P \cap Q_\la - \{(v,0)\} \eqd \la^{1/(d-1)}([ \tP_\la \cap Q_1 ] - \{(v',0)\}),
  $$
  where $v = \la^{1/(d-1)} v'$.   We thus obtain by translation invariance of $\xi$
  $$
\la^{-1} \E [\langle 1, \mu_\la \rangle] = \int_{\tQ_1} \int_{- \infty}^{\infty}
 \E \xi( (\0,h_0), \la^{1/(d-1)}([ \tP_\la \cap Q_1 ] - \{(v_0',0)\})  e^{h_0} dh_0 dv_0'.
 $$




Notice that for all $v_0' \in \tQ_1$ we have $\la^{1/(d-1)}([ \tP_\la \cap Q_1 ] - \{(v_0',0)\}) \tod \P$ as $\la \to \infty$. Also,
the functional $\xi$ satisfies the spatial localization and moment conditions
of those functionals in $\Xi^{(\infty)}$ and consequently we have for all $v_0' \in \tQ_1$
$$
\lim_{\la \to \infty} \E \xi( (\0,h_0), \la^{1/(d-1)}([ \tP_\la \cap Q_1 ] - \{(v_0',0)\})) =
\E [\xi( (\0,h_0), \P)].
$$
As in  Lemma \ref{L1}, we may show that $\E [\xi( (\0,h_0), \la^{1/(d-1)}([ \tP_\la \cap Q_1 ] - \{(v_0',0)\}))] e^{h_0}$ is dominated by
an exponentially decaying function of $h_0$, uniformly in $\la$. Thus by the dominated convergence theorem we get
\be \label{LimitY}
\lim_{\la \to \infty} \la^{-1} \E [ \langle 1, \mu_\la \rangle ] =
 \int_{- \infty}^{\infty}
\E [\xi( (\0,h_0), \P))] e^{h_0} dh_0.
\ee

To prove variance asymptotics, we argue as follows. 
For all $h_0 \in \R$,  and $(v_1,h_1) \in \R^{d-1} \times \R,$
we abuse notation and put
$$
c^{\xi}((\0,h_0), (v_1,h_1), \P\cap Q_\la ) := \E[ \xi((\0,h_0), (\P\cap Q_\la) \cup \{(v_1,h_1)\}) \times \xi((v_1,h_1),(\P\cap Q_\la) \cup \{(\0, h_0)\})]
 $$
$$- \E [(\xi((\0,h_0), \P\cap Q_\la)] \E
[\xi((v_1,h_1), \P\cap Q_\la )].
$$
Then we have
\begin{align}
\hspace*{-1cm}\ \  \ \ & \ \ \la^{-1}\Var [\langle 1, \mu_\la \rangle] \nonumber\\
 & = \la^{-1} \int_{\tQ_\la} \int_{- \infty}^{\infty}
 \E \left[\xi( (v,h), P\cap Q_\la )\right] e^h dh_0 dv_0 \nonumber\\
  & \ \ + \la^{-1} \int_{\tQ_\la} \int_{- \infty}^{\infty} \int_{\tQ_\la} \int_{- \infty}^{\infty}
 c^{\xi}((v_0,h_0), (v_1,h_1), P\cap Q_\la) e^{h_0 + h_1} dh_0 dh_1 dv_1 dv_0 \nonumber\\
& =  \int_{\tQ_1}\int_{- \infty}^{\infty} \E \left[\xi( (\0,h_0), \la^{1/(d-1)}([ \tP_\la \cap Q_1 ] - \{(v_0',0)\}))\right]  e^{h_0} dh_0 dv_0'
\nonumber \\
& \ \ +  \int_{\tQ_1}\int_{- \infty}^{\infty} \int_{\tQ_\la - \lambda^{1/(d-1)}v_0'} \int_{- \infty}^{\infty}
 c^{\xi}((\0,h_0), (v_1,h_1), \la^{1/(d-1)}([ \tP_\la \cap Q_1 ] - \{(v_0',0)\}))\nonumber\\&\hspace*{12cm} e^{h_0 + h_1} dh_0 dh_1 dv_1 dv_0'.\label{corovar}
\end{align}

The first integral  in \eqref{corovar} converges to the limit \eqref{LimitY}.
We show that the second term in \eqref{corovar} converges.  For all $h_0, h_1, v_1$ and $v_0'$ we have
$$
\lim_{\la \to \infty} c^{\xi}((\0,h_0), (v_1,h_1), \la^{1/(d-1)}([ \tP_\la \cap Q_1 ] - \{(v_0',0)\})) e^{h_0 + h_1} = c^{\xi}((\0,h_0), (v_1,h_1), \P)e^{h_0 + h_1} ,
$$
where  $c^{\xi}$ is defined as in \eqref{SO2}.
As in the proof of Lemma \ref{L3}, we have  that
$$ c^{\xi}((\0,h_0), (v_1,h_1), \la^{1/(d-1)}([ \tP_\la \cap Q_1 ] - \{(v_0',0)\})) e^{h_0 + h_1}$$
is dominated by an exponentially decaying function of $h_0$, $h_1$, and $v_1$, uniformly in $\la \in [\la_0, \infty)$.
Consequently, by the dominated convergence theorem, we obtain
$$
\liml \la^{-1} \Var [\langle 1, \mu_\la \rangle] =  \int_{- \infty}^{\infty}
\E [\xi( (\0,h_0), \P))] e^{h_0} dh_0 \ +
$$
$$
+ \int_{- \infty}^{\infty} \int_{\R^{d-1}} \int_{- \infty}^{\infty}
c^{\xi}((\0,h_0), (v_1,h_1), \P)e^{h_0 + h_1}dh_0 dv_1 dh_1.
$$
This concludes the proof of Corollary  \ref{Th3}.  \qed


\noindent{\bf 5.5. Proof of Theorem \ref{intrin}.}
Let us fix $k\in \{1,\cdots,d-1\}$. We follow the same method and notation as on pages 54-55 in \cite{CSY}. The key idea is to use Kubota's formula (see (6.11) in \cite{SW}) which says roughly that the $k$-th intrinsic volume of $K_{\la}$ is,  up to a multiplicative constant,  equal to the mean over the set $G(d,k)$ of all $k$-dimensional linear subspaces of $\R^d$ of the $k$-dimensional Lebesgue measure of the projection of $K_{\la}$ onto $L$, denoted by $K_{\la}|L$. In other words, we have
        \begin{equation}\label{KUBOTA}
         V_k(K_{\la}) = 
         \frac{d! \kappa_d}{k! \kappa_k (d-k)! \kappa_{d-k}}
         \int_{G(d,k)} \Vol_k(K_{\la} | L) d\nu_k(L),
        \end{equation}
        where $\nu_k$ is the normalized Haar measure on the $k$-th Grassmannian  $G(d,k)$ of $\R^d$. 

For every $x\in\R^d\setminus\{\0\}$, We consider now the function
$\vartheta_{L}(x,K_{\la})={\bf 1}(x\not\in K_{\la}|L)$ and the so-called projection avoidance functional
$$\vartheta_k(x,K_{\la})=\int_{G(\mbox{\tiny{lin}}[x],k)}\vartheta_{L}(x,K_{\la})d\nu_k^{\mbox{\tiny{lin}}[x]}(L)$$
where
$\mbox{lin}[x]$ is the one-dimensional linear space spanned by $x$, $G(\mbox{lin}[x],k)$ is the set of $k$-dimensional linear subspaces of $\R^d$ containing $\mbox{lin}[x]$ and $\nu_k^{\mbox{\tiny{lin}}[x]}$ is the normalized Haar measure on $G(\mbox{lin}[x],k)$ (see (2.7) in \cite{CSY}).
Using \eqref{KUBOTA} and Fubini's theorem, we rewrite the defect intrinsic volume of $K_{\la}$ as
\begin{align*} \label{eq:intrinsicvolumeformula}
 V_k(B_d(\0,R_{\la}))-V_k(K_{\la})&=
\frac{\binom{d-1}{k-1}}{\kappa_{d-k}}\int_{\R^d}[\vartheta_k(x,K_{\la})
-\vartheta_k(x,B_d(\0,R_{\la}))]\frac{dx}{|x|^{d-k}}.
\end{align*}
In particular, we have the decomposition
$$ V_k(B_d(\0,R_{\la}))-V_k(K_{\la})=R_\la^{-(d+1-k)}\sum_{x\in \P_{\la}} \xi_{V,k}(x,\P_{\la})$$
where
 \begin{equation}
   \label{eq:xi_vk}
 \xi_{V,k}(x,\P_{\la}):=d^{-1}\frac{\binom{d-1}{k-1}}{\kappa_{d-k}}
 \int_{{\rm{cone}}(x, \P_\la) }
 [\vartheta_k(y,K_{\la})-
 \vartheta_k(y,B_d(\0,R_{\la}))]\frac{R_\la^{d+1-k}dy}{|y|^{d-k}}
 \end{equation}
if $x$ is extreme and $\xi_{V,k}(x,\P_{\la})=0$ otherwise. The corresponding empirical measure is
$$\mu^{\xi_{V,k}}_\la:= \sum_{x \in \P_\la} \xi_{V,k}(x, \P_\la) \de_x.$$
In the same spirit as in the definition of the defect volume functional (see Definition \ref{vol}), the scaling factor $R_{\la}^{d+1-k}$ is artificially inserted inside the integral in \eqref{eq:xi_vk} in order to use later the convergence \eqref{eq:convmk}. In view of the equality
$$\langle 1, \mu_\la^{\xi_{V,k}} \rangle=R_{\la}^{d+1-k}[V_k(B_d(\0,R_{\la}))-V_k(K_{\la})],$$
we observe that it is enough to show that \eqref{main1} and \eqref{main2} are satisfied when $\xi=\xi_{V,k}$.

We notice that the equalities \eqref{dis1} and \eqref{dis3} hold when $\xi$ is put to be  $\xi_{V,k}$. Let us define $$\xi_{V,k}^{(\la)}(w,\P^{(\la)})=\xi_{V,k}([T^{(\la)}]^{-1}(w), \P_{\la}),\quad w\in \R^{d}.$$ Observe from the proof of Theorem \ref{Th5} that it is enough to show the convergence up to a multiplicative rescaling of each of the terms $\E[\xi_{V,k}^{(\la)}(w,\P^{(\la)})]$, $\E[\xi_{V,k}^{(\la)}(w,\P^{(\la)})^2]$ and $c^{\xi_{V,k}^{(\la)}}(w,w')$ where $w,w'\in\R^d$, as well as bounds similar to those in Lemmas \ref{L1} and \ref{L3}.

Let us show for instance the convergence of $\E\xi_{V,k}^{(\la)}(w,\P^{(\la)}), w \in W_\la$. The localization radius associated with $\xi_{V,k}^{(\la)}$ is the same as that for $\xi_V^{(\la)}$. Moreover, a moment bound similar to \eqref{LIMITBD2} can be obtained when $\xi^{(\lambda)}$ is replaced by $\xi_{V,k}^{(\lambda)}$. We need to introduce now the scaling limit $\xi_{V,k}^{(\infty)}(w,\P)$ of $\xi_{V,k}^{(\la)}(w,\P^{(\la)})$, i.e. the exact analogue of \eqref{eq:xi_vk} when Euclidean convex geometry is replaced by parabolic convex geometry (see e.g. section 3 in \cite{CSY}). For every $w=(v,h)\in\R^{d-1}\times \R$, we denote by $w^{\updownarrow}$ the set $\{v\}\times \R$ and by $\mu^{w^{\updownarrow}}_k$ the normalized Haar measure on the set $A(w^{\updownarrow},k)$ of all $k$-dimensional affine spaces in $\R^d$ containing $w^{\updownarrow}$. For every affine space $L$ containing $w^{\updownarrow}$, we define the corresponding orthogonal paraboloid volume  $\Pi^{\perp}[w;L]$ as the set $\{w'=(v',h')\in \R^{d-1}\times\R:(w-w')\perp L,h'\le h-\frac{\|v-v'\|^2}{2}\}$. In other words, $\Pi^{\perp}[w;L]$ is the set of points of $(w\oplus L^{\perp})$ which are under the paraboloid surface $\partial \Pi^{\downarrow}(w)$ with apex at $w$.  We put $\vartheta_L^{(\infty)}(w)=1$ when $\Pi^{\perp}[w;L]\cap T^{(\lambda)}(K_{\lambda})= \emptyset$ and $0$ otherwise. In particular, when $L=w^{\updownarrow}$, we have $\vartheta_{w^{\updownarrow}}^{(\infty)}(w)={\bf 1}(\Pi^{\perp}[w;L]\cap \Phi=\emptyset).$ Now for every $w\in\R^d$, we define $$\vartheta_k^{(\infty)}(w):=\int_{A(w^{\updownarrow},k)}\vartheta_L^{(\infty)}(w)d\mu^{w^{\updownarrow}}_k(L).$$ We finally define the limit score $$\xi_{V,k}^{(\infty)}(w,\P):=d^{-1}\frac{\binom{d-1}{k-1}}{\kappa_{d-k}}\int_{{\rm{Cyl}}(w)\times \R}[\vartheta_k^{(\infty)}(w')-{\bf 1}(\{w'\in\R^{d-1}\times\R_-\})]dw'$$
if $w\in \mbox{Ext}(\P)$ and $\xi_{V,k}^{(\infty)}(w,\P) =0$ otherwise.
Denote by $m_k$ the measure $|x|^{k-d}dx$. In the same spirit as Lemma \ref{weakcon}, have as $\la \to \infty$
\begin{equation}
  \label{eq:convmk}
T^{(\la)}(R_{\la}^{d+1-k}dm_k) \tod \Vol_d.
\end{equation}
As in the proof of Lemma \ref{L2a}, we may show via the continuous mapping theorem that for every $w\in\R^d$, we have the convergence in distribution $\xi_{V,k}^{(\la)}(w,\P^{(\la)})\overset{D}{\to} \xi_{V,k}^{(\infty)}(w,\P)$. We deduce from the moment bound the analog of Lemma \ref{L2a} for  $\xi_{V,k}^{(\la)}$. In the spirit of Lemma \ref{L2}, we finally
show that
$\lim_{\la\to\infty}\E\xi_{V,k}^{(\lambda)}(w,\P^{(\la)})
=\E\xi_{V,k}^{(\infty)}(w,\P).$

\vskip.5cm

\noindent{\em Acknowledgements}. 
We thank an anonymous referee for remarks leading to an improved exposition. P. Calka also 
thanks the Department of Mathematics at Lehigh University for its kind hospitality and support.

\noindent Pierre Calka, Laboratoire de Math\'ematiques Rapha\"el Salem, Universit\'e de Rouen,
 Avenue de l'Universit\'e, BP.12, Technop\^ole du Madrillet, F76801 Saint-Etienne-du-Rouvray France;
\ \ \  {\texttt pierre.calka@univ-rouen.fr}

\vskip.5cm

\noindent J. E. Yukich, Department of Mathematics, Lehigh University,
Bethlehem PA 18015, USA; \ \ {\texttt joseph.yukich@lehigh.edu}

\end{document}